\documentstyle[12pt,amsfonts,fullpage]{article}

\def\SBIMSMark#1#2#3{
 \font\SBF=cmss10 at 10 true pt
 \font\SBI=cmssi10 at 10 true pt
 \setbox0=\hbox{\SBF Stony Brook IMS Preprint \##1}
 \setbox2=\hbox to \wd0{\hfil \SBI #2}
 \setbox4=\hbox to \wd0{\hfil \SBI #3}
 \setbox6=\hbox to \wd0{\hss
             \vbox{\hsize=\wd0 \parskip=0pt \baselineskip=10 true pt
                   \copy0 \break%
                   \copy2 \break%
                   \copy4 \break}}
 \dimen0=\ht6   \advance\dimen0 by \vsize \advance\dimen0 by 8 true pt
                \advance\dimen0 by -\pagetotal
 \dimen2=\hsize \advance\dimen2 by .25 true in
 \ht6=0pt \dp6=0pt
 \setbox8=\vbox to \dimen0{\vfill \hbox to \dimen2{\hss \copy6}}
 \ht8=0pt \dp8=0pt \wd8=0pt
 \copy8
 \message{*** Stony Brook IMS Preprint #1, #2 ***}
}

\begin{document}



\newcommand{\R}{{\Bbb R}}
\newcommand{\ct}{\centerline}

\newcommand{\ga}{\alpha}
\newcommand{\gga}{\gamma}
\newcommand{\gG}{\Gamma}
\newcommand{\gb}{\beta}
\newcommand{\gd}{\delta}
\newcommand{\gD}{\Delta}
\newcommand{\gk}{\kappa}
\newcommand{\get}{\eta}
\newcommand{\gep}{\varepsilon}
\newcommand{\gvp}{\varphi}
\newcommand{\gf}{\varphi}
\newcommand{\gl}{\lambda}
\newcommand{\gL}{\Lambda}
\newcommand{\Gl}{\Lambda}
\newcommand{\gch}{\chi}
\newcommand{\gp}{\pi}
\newcommand{\gps}{\psi}
\newcommand{\gs}{\sigma}
\newcommand{\gS}{\Sigma}
\newcommand{\gt}{\theta}
\newcommand{\Gt}{\Theta}
\newcommand{\gT}{\Theta}
\newcommand{\gw}{\omega}
\newcommand{\gW}{\Omega}
\newcommand{\gx}{\xi}
\newcommand{\gz}{\zeta}

\newcommand{\bA}{\hbox{\bf A}}
\newcommand{\ba}{\hbox{\bf a}}
\newcommand{\bB}{\hbox{\bf B}}
\newcommand{\bb}{\hbox{\bf b}}
\newcommand{\bC}{\hbox{\bf C}}
\newcommand{\bc}{\hbox{\bf c}}
\newcommand{\bD}{\hbox{\bf D}}
\newcommand{\bE}{\hbox{\bf E}}
\newcommand{\bff}{\hbox{\bf f}} 
\newcommand{\bH}{\hbox{\bf H}}
\newcommand{\bh}{\hbox{\bf h}}
\newcommand{\bI}{\hbox{\bf I}}
\newcommand{\bk}{\hbox{\bf k}}
\newcommand{\bL}{\hbox{\bf L}}
\newcommand{\bl}{\hbox{\bf l}}
\newcommand{\bM}{\hbox{\bf M}}
\newcommand{\bN}{\hbox{\bf N}}
\newcommand{\bP}{\hbox{\bf P}}
\newcommand{\bp}{\hbox{\bf p}}
\newcommand{\bQ}{\hbox{\bf Q}}
\newcommand{\bq}{\hbox{\bf q}}
\newcommand{\bR}{\hbox{\bf R}}
\newcommand{\br}{\hbox{\bf r}}
\newcommand{\bS}{\hbox{\bf S}}
\newcommand{\bT}{\hbox{\bf T}}
\newcommand{\bt}{\hbox{\bf t}}
\newcommand{\bU}{\hbox{\bf U}}
\newcommand{\bu}{\hbox{\bf u}}
\newcommand{\bv}{\hbox{\bf v}}
\newcommand{\bx}{\hbox{\bf x}}
\newcommand{\by}{\hbox{\bf y}}
\newcommand{\bZ}{\hbox{\bf Z}}
\newcommand{\bz}{\hbox{\bf z}}

\newcommand{\CA}{{\cal A}}
\newcommand{\CB}{{\cal B}}
\newcommand{\CC}{{\cal C}}
\newcommand{\CD}{{\cal D}}
\newcommand{\CE}{{\cal E}}
\newcommand{\CF}{{\cal F}}
\newcommand{\CG}{{\cal G}}
\newcommand{\CH}{{\cal H}}
\newcommand{\CI}{{\cal I}}
\newcommand{\CJ}{{\cal J}}
\newcommand{\CK}{{\cal K}}
\newcommand{\CL}{{\cal L}}
\newcommand{\CM}{{\cal M}}
\newcommand{\Cm}{{\cal m}}
\newcommand{\CN}{{\cal N}}
\newcommand{\CO}{{\cal O}}
\newcommand{\CP}{{\cal P}}
\newcommand{\CQ}{{\cal Q}}
\newcommand{\CR}{{\cal R}}
\newcommand{\CS}{{\cal S}}
\newcommand{\CT}{{\cal T}}
\newcommand{\CU}{{\cal U}}
\newcommand{\CV}{{\cal V}}
\newcommand{\CW}{{\cal W}}
\newcommand{\CX}{{\cal X}}
\newcommand{\CY}{{\cal Y}}
\newcommand{\CZ}{{\cal Z}}

\newcommand{\I}{{\bf I}}
\newcommand{\II}{{\bf I\kern -1pt I}}

\newcommand{\A}{\forall}
\newcommand{\E}{\exists}

\newcommand{\titles}[1]{\bigskip\ct {\Large\bf #1}\bigskip \bigskip}
\newcommand{\titless}[1]{\bigskip\ct {\Large\bf #1}\bigskip }
\newcommand{\sections}[1]{\addtocounter{section}{1}\setcounter{theorem}{0}
\bigskip \nin {\large\bf #1}\medskip\medskip}
\newcommand{\authors}[3]{\bigskip\ct {#1} \bigskip\ct{#2}
\bigskip\ct{#3} \bigskip}
\newcommand{\authorss}[4]{\bigskip\ct {#1} \ct{#2}
\ct{#3} \ct{#4} }

\SBIMSMark{1996/9}{August 1996}{}

\vspace{3cm}
\titless{A new advance in the Bernstein Problem}
\titless{ in mathematical genetics}

\bigskip

\ct{\rm Yuri I. Lyubich}

\ct{\rm Department of Mathematics, Technion}
\ct{32000 Haifa, Israel}
\bigskip


\medskip
\section{Introduction}

$\ \ \ \ \ \ \ \ \ $Here we present a new result on a problem 
posed by S.N.Bernstein [1]  in mathematical 
foundations of the population genetics. The problem is related to a 
statement called 
{\em The Stationarity Principle (S.P.)}. Being valid under the Mendel Law 
this principle is consistent with some more general mechanisms of heredity.
Bernstein suggested to describe all situations satisfying the S.P.. 
In mathematical terms this
sounds as follows.

Let $\:\gD^{n-1}\:\subset\:\R^n\:$ be the basis simplex,
$$
\gD^{n-1}\:=\:\{\:x\:\in\:\R^n:\ s(x)\:\equiv\:
\sum\limits_{i=1}^nx_i=1,\ \ x\:=\:(x_i)^n_1\:\geq\:0\}.
\eqno (1.1)
$$

Consider a quadratic mapping $\:V\: :\:\gD^{n-1}\:\rightarrow\:\gD^{n-1},\:$
$$
x^{'}_j\:\equiv\:(Vx)_j\:=\:\sum\limits_{i,k=1}^np_{ik,j}x_ix_k\ \ \ 
(1\:\leq\:j\:\leq\:n),
\eqno (1.2)
$$
which is {\em stochastic} in the sense that
$$
p_{ik,j}\:\geq\:0,\ \ \ \sum\limits_{j=1}^np_{ik,j}\:=\:1.
\eqno (1.3)
$$
Certainly, the symmetry $\:p_{ki,j}\:=\:p_{ik,j}\:$ is also supposed.

A mapping
$\:V\:$ is called {\em Bernstein} (or {\em stationary}) if 
$\:V^2\:=\:V,\:$ where $V^2\equiv V\circ V$.
This property is just the S.P. 
{\em The Bernstein problem is to explicitly
describe all such mappings.} 
For $\:n\:=\:3\:$ the problem was solved  
in [2], [3], [4]. The cases $\:n\:=\:1,2\:$ are trivial (see below).

Biologically, $\:V\:$ is the  evolutionary operator of an infinite
 population under certain conditions (see [17], Sections 1.1,
1.2 for a detailed  explanation ).
 Each individual from the population belongs to a 
biological {\em type (character)}. The set of types is supposed to be 
finite, say $\:\{\:1,\: ...,\:n\:\},\:$ and the partition of 
the population into the types has to be {\em hereditary}. This means that
for every triple $\:(i,k,j)\:$ of types there exists a probability 
$\:p_{ik,j}\:$ for parents of types $\:i\:$ and $\:k\:$ 
to have an offspring of type $\:j\:$. In this sense 
$\:p_{ik,j}\:$ are the {\em inheritance coefficients.} Thus, we have (1.3)
automatically in this context. The symmetry $\:p_{ki,j}\:=\:p_{ik,j}\:$
means that the sexual differentiation does not affect on the heredity.

The points $\:x\:\in\:\gD^{n-1}\:$ are just the probability distributions 
on the set of types. Every such a point is a {\em state} of the population.
If $\:x\:$ is a state in a {\em parental generation} then $\:x^{'}\:=\:Vx\:$
is the state in the {\em offspring generation}. With an {\em initial state} 
$\:x\:$
the sequence $\:\{\:V^tx\:\}^{\infty}_{t=0}\:$ is the corresponding {\em
trajectory} of the population considered as a dynamical system ([17],
Section 1.2). In a simplest case $\:x\:$ is a fixed point, $\:x=Vx,\:$ so
the trajectory is reduced to the point $\:x.\:$ Such points are {\em
equilibria} from the dynamical point of view. 
The only case all states are equilibria
is $\:V\:=\:I,\:$ the identity mapping. Note that this mapping $\:x^{'}_j\:=
\:x_j\ \ (1\:\leq\:j\:\leq\:n)\:$ can be also represented as a quadratic one,
namely, $\:x^{'}_j\:=\:x_js(x)\ \ (1\:\leq\:j\:\leq\:n).\:$ This is also
can be done for any linear stochastic mapping $\:x^{'}\:=\:Tx,\:$ namely,
$\:x^{'}\:=\:s(x)Tx.\:$ The corresponding dynamical system is the Markov
chain generated by $\:T.\:$ If $\:(t_{ij})^n_{i,j=1}\:$ is the matrix of
$\:T\:$ then $\:p_{ik,j}\:=\:\frac{1}{2}(t_{is}\:+\:t_{kj})\:$ in the 
above mentioned quadratic representation.

It is also useful to note that any constant mapping 
$\:x^{'}\:=\:c\:\in\:\gD^{n-1}
\:$ can be written as a quadratic one: $\:x^{'}\:=\:cs^2(x),\ \ 
x\:\in\:\gD^{n-1}.$

It is easy to prove that {\em for $\:n\leq 2\:$ every Bernstein mapping is
constant or identity}.

The S.P. $\:V^2\:=\:V\:$ means that every offspring state $\:Vx\:$ is an 
equilibrium, so that the trajectory consists of $\:x\:$ and $\:Vx.\:$
Such a simplest dynamics should correspond to an ``elementary" law of 
heredity.
Just this philosophy compeled S.N.Bernstein to pose his problem. On
the other hand, S.P. is a fortiori valid under the Mendel Law. Let us 
explain this in a more detail.

A simplest mechanism of heredity (discovered by Mendel) is determined by 
two {\em genes}, say $\:{\bf A}\:$ and $\:{\bf a}.\:$ Every individual has 
one of three possible genotypes: $\:{\bf AA},\:{\bf aa},\:{\bf Aa}.\:$ 
Each  parent
provides each offspring with one of these two  genes. 
The genes of $\:{\bf Aa}\:$ are reproduced in offspring with  
probabilities $\:\frac{1}{2}.\:$ Any offspring genotype appears as an
independent random  combination of two parental genes. This mechanism
transforms a parental state $\:x\:=\:(x_1,x_2,x_3)\:$ of a population 
into the offspring state $\:x^{'}\:$ with
$$
x^{'}_1\:=\:p^2,\ \ x^{'}_2\:=\:q^2,\ \ x^{'}_3\:=\:2pq
\eqno (1.4)
$$
where
$$
p\:=\:x_1\:+\:\frac{1}{2}x_3,\ \ q\:=\:x_2\:+\:\frac{1}{2}x_3
\eqno (1.5)
$$
These $\:p\ {\rm and}\ q\:$ are the probabilities of the genes 
$\:{\bf A}\:$ and $\:{\bf a}\:$ at the state 
$\:x.\:$ For the first time these formulas were independently obtained
in [5] and [18] therefore the corresponding quadratic mapping $\:V\: :\:
\gD^2\:\rightarrow\:\gD^2\:$ is called the {\em Hardy-Weinberg mapping.}
Obviously, $\:p+q=1,\:$ so in the next generation
$$
p^{'}\:=\:x^{'}_1\:+\:\frac{1}{2}x^{'}_3\:=\:p^2\:+\:pq\:=\:p
\eqno (1.6)
$$
and $\:q^{'}\:=\:q\:$ similarly. For this reason 
$\:x^{''}_1\:=\:p^{{'}^2}\:=\:p^2\:=\:x^{'}_1\:$ and 
$\:x^{''}_2\:=\:x^{'}_2,\ \ x^{''}_3\:=\:x^{'}_3.\:$ This means that  
$\:V^2\:=\:V.\:$

The relations $\:p^{'}=p\:$ and $\:q^{'}=q\:$ catch the phenomenon that  the 
genes pass from parents to offsprings with no appearing or disappearing 
(usually occured under mutation and selection).

Following this classical pattern we introduced ([7]; [17],
Section 4.1) a general concept of {\em stationary gene structure.} Let 
us reproduce it below.

Given an evolutionary operator $\:V\: :\:\gD^{n-1}\: 
\rightarrow\:\gD^{n-1},\:$ a linear
form
$$
f(x)\:=\:\sum\limits_{i=1}^na_ix_i
\eqno (1.7)
$$
is called {\em invariant} if $\:f(Vx)\:=\:f(x)\ \ (x\:\in\:\gD^{n-1})\:$
or $\:f^{'}=f\:$ for short, like (1.6). A trivial example is $\:s(x)\:$ 
or any multiple of  $\:s\:$. Obviously, the set $\:J\:$ of all invariant
linear forms is a linear space, 
$\:1\leq {\rm dim}J\leq n.\:$ If $\:V\:$ is constant then 
dim$J=1.\:$ If $\:V=I\:$ (and only in this case) then dim$J=n.\:$

We say that $\:V\:$ {\em has a stationary gene structure (s.g.s)} 
($\:\equiv\:V\:$ is {\em regular}) if this
mapping can be written as
$$
x^{'}_j\:=\:\sum\limits_{i,k=1}^rc_{ik,j}f_i(x)f_k(x)\ \ \ 
(1\:\leq\:j\:\leq\:n)
\eqno (1.8)
$$
where $\:f_1,...,f_r\:$ are some invariant linear forms. 
If so, these forms can be chosen in a special {\em canonical} way. 
Namely, one can consider the cone $\:C\:$  of those $\:f\:\in\:J\:$
which are {\em nonnegative} in the usual sense: all $\:a_i\geq 0\:$ in
(1.7). This is a closed polyhedral cone; $\:s\:\in\:{\rm Int}C,\:$ so
$\:C\;$ generates the space $\:J.\:$ In (1.8) $\:f_1, ..., f_r\:$ can
be taken from the extremal rays of $\:C,\:$ one for each ray. In such
a way let
$$
f_l(x)\:=\:\sum\limits_{j=1}^n\pi_{lj}x_j\ \ \ (1\:\leq\:l\:\leq\:r)
\eqno (1.9)
$$
with
$$
{\rm max}_j\pi_{lj}\:=\:1\ \ \ (1\:\leq\:l\:\leq\:r)
\eqno (1.10)
$$
Then we say that $\:f_l\:$ are {\em canonical.}

Obviously, $\:r\:\geq\:{\rm dim}J\:$ and $\:r\:=\:{\rm dim}J\:$ iff
$\:f_1, ..., f_r\:$ are linearly independent.
In the last case we say that $\:V\:$ has an {\em elementary gene 
structure (e.g.s).} Then the coefficients $\:c_{ik,j}\:$ in (1.8) are 
uniquely determined and $\:c_{ik,j}\:\geq\:0\:$ ([7]; [17], Corollary 4.3.4
and Lemma 4.3.5). For {\em nonelementary gene structure (n.e.g.s.) 
$\:c_{ik,j}\:$ are not uniquely determined but there exists a set of
$\:c_{ik,j}\:\geq\:0\:$ in (1.8) with canonical $\:f_1,\: ...,\:f_r\:$ }
([12]; [17], Theorem 4.4.1).

A simplest example of n.e.g.s came from [4]. This is the {\em quadrille 
mapping} $\:V\: :\:\gD^3\:\rightarrow\:\gD^3,\:$ 
$$
x^{'}_1\:=\:p_1q_1,\ \ x^{'}_2\:=\:p_2q_2,\ \ x^{'}_3\:=\:p_1q_2,\ \ 
x^{'}_4\:=\:p_2q_1
\eqno (1.11)
$$
with
$$
p_1\:=\:x_1+x_3,\ \ p_2\:=\:x_2+x_4,\ \ 
q_1\:=\:x_1+x_4,\ \ q_2\:=\:x_2+x_3.
\eqno (1.12)
$$
A relevant genetical mechanism was suggested in [7].

Like the Hardy-Weinberg case we have $\:V^2\:=\:V\:$ {\em for any s.g.s.}
The converse is not true ([17], p.172).

A {\em restricted Bernstein problem} posed in [7] {\em is to 
explicitly describle all regular stochastic quadratic mappings.} 
In the case of e.g.s. this
problem was solved in [7]. Later the general case was 
solved  in [8], [11], [12] (see also [17], Chapter 4). However, it turned 
out that a satisfactory genetical interpretation reguires an additional
property of {\em normality}.

A stochastic quadratic mapping $\:V\: 
:\:\gD^{n-1}\:\rightarrow\:\gD^{n-1}\:$ is called {\em normal} if in (1.2)

1) all $\:x^{'}_j\:\not\equiv\:0\:$ ({\em the nondegeneracy)};

2) every two $\:x^{'}_{j_1},x^{'}_{j_2}\ \ (j_1\neq j_2)\:$ are not 
proportional ({\em the external irreducibility)};

3) there is no pair $\:i,k\ \ (i\neq k)\:$ such that all $\:x^{'}_j\:$ only
depend on $\:x_i+x_k\:$ and $\:x_l\ \ (l\neq i,\ l\neq k)\:$ ({\em the
internal irreducibility}).

A constant mapping is normal in the only case $\:n=1.\:$ The unit mapping 
is normal in all dimensions $\:n.\:$  

If $\:V\:$ is not normal one can reduce it to a normal one by a standard
procedure of {\em normalization} ([11]; [17], Section 3.9). Under the 
normalization the dimension $\:n\:$ decreases but this process
preserves the regularity, moreover, $\:r\:$ and dim$J$ are invariant.

The explicit description of all regular normal evolutionary operators 
is contained in [17], Theorem 4.3.9 (for e.g.s.) and Theorem 4.6.1 (for 
n. e.g.s.). We also 
explain this below (Section 3) in a more apparent algebraic form remarkably 
corresponding to some genetical mechanisms (cf. [17], p.189,207).
The point is that the types $\:\{1,\: ....\:n\}\:$ in any normal s.g.s can
be identified with
some pairs of genes $\:{\bf A}_1,\:,,,\: {\bf A}_r;\:$ the probabilities 
of these genotypes at a state $\:x\:$ are the canonical 
$\:f_l(x).\:$ If a s.g.s. is 
not normal (but nondegenerate) then there are some different types 
whose formal genotypes are the same, so some of the types are redundant 
([17], Section 4.2).

The degeneracy means that some of types disappear from the population 
after mating. Such types can not be considered as hereditarily significant 
ones.

Thus, in the Bernstein problem the only case of normal s.g.s. has a 
genetical sense. 
However, in this context the S.P. is not an axiom, it is 
a consequence of s.g.s.
If we wish to preserve the S.P. as an axiom then we should add 
something else to get s.g.s. as a consequence providing a natural genetical
interpretation. In such a way some relevant conjectures were suggested in 
[16] (see also [17], Section 5.7). A proof of one of them is  the subject
of the present paper.

{\bf Definition.} {\em A stochastic quadratic mapping $\:V\:$ is called
ultranormal if its restrictions to all invariant faces of the simplex 
$\:\gD^{n-1}\:$ are normal.}

The ultranormality is a natural axiom in addition to the S.P. since we 
recognize the normality as a necessary property of {\em all} evolutionary 
operators, in particular, of the restrictions of $\:V\:$ to all 
invariant faces.

{\bf Main Theorem.} {\em Every ultranormal stochastic Bernstein mapping
$\:V\:$ is regular.}

Our above mentioned works contain a proof of this theorem in the case 
dim(Im$V)\:\leq\:2\:$ or $\geq\:n-2,\:$ in particular, for $n\:\leq\:5.\:$

The Main Theorem combining with our explicit description 
of the regular normal mappings completely resolves the Bernstein problem in 
the ultranormal case. Note that {\em every normal s.g.s. is ultranormal}
as directly follows from its explicit form.

Is  every {\em normal} stochastic Bernstein mapping regular? This is an 
open question for $\:n\:\geq\:5\:$ (cf. [17], Section 5.7).
An  affirmative answer would be a key to the general Bernstein problem by 
normalization.

Our approach to the Main Theorem is basically algebraical and partly 
topological one. In Section 
3.4 some relevant means are prepared. The corresponding key words are 
``Bernstein algebra",``regular algebras", 
``stochastic algebras" and  their ``offspring 
subalgebras".  

The proof of the Main Theorem is given in Section 5. Actually, 
we prove that {\em every ultranormal stochastic Bernstein algebra is regular}
or, equivalently,  admits the above mentioned explicit form.

The result of this paper  was announced  at the 9th Haifa Matrix 
Theory Conference   on June  1, 1995.

\smallskip\noindent{\bf Aknowledgemnt.}
This paper was partly prepared during the author's visit
to the IMS at Stony Brook in the summer 1995. 
The author thanks the IMS for their  hospitality

\section{Bernstein algebras.}
$\ \ \ \ \ \ \ \ \ $
For any evolutionary mapping $\:V\:$ one can consider an algebra
$\:\CA_V\:$ in $\:\R^n\:$ whose structure constants at the canonical basis 
$\:\{\:e_j\:\}^n_1\:$ are $\:p_{ik,j},\:$ so that we have the 
multiplicative table
$$
e_ie_k\:=\:\sum\limits_{j=1}^np_{ik,j}e_j.
\eqno (2.1)
$$

The algebra $\:\CA_V\:$ is commutative but, as rule, it is not associative.

In a biological interpretation, the types $\:\{1,\:...,\:n\:\}\:$ have 
to be identified with the corresponding basis vectors $\:\{\:e_1,\:...,\:
e_n\:\}.\:$ With parental types $\:e_i,e_k\:$ an offspring is of  
type $\:e_j\:$ with probability $\:p_{ik,j}.\:$   

The {\em evolutionary algebra} is {\em stochastic} in the sense that the 
symplex $\:\gD^{n-1}\:$ is invariant with respect
to the multiplication. Indeed, $\:x\:\geq\:0\:\& \:y\:\geq\:0\:\Rightarrow\:
xy\:\geq\:0\:$ 
and $\:s(x)\:=\:1\:\&\\ 
\:s(y)\:=\:1\:\Rightarrow\:s(xy)\:=\:1\:$
because $\:s\:$ is a multiplicative linear functional (a {\em weight)}:
$$
s(xy)\:=\:s(x)s(y).
\eqno (2.2)
$$
This means that the pair $\:(\CA_V,s)\:$ is a real {\em baric} algebra (see
[14]; [17],Sections 3.3, 3.8).

Note that $\:Vx\:=\:x^2 \ \ (x\:\in\:\gD^{n-1})\:$ so $\:x^2\:$ is a (unique)
quadratic extension $\:\tilde{V}\:$ of $\:V\:$ from $\:\gD^{n-1}\:$ to the 
whole 
space $\:\R^n.\:$ It is very fruitfull to reformulate the Bernstein 
problem algebraically. We systematically used this approach 
earlier starting with the following

{\bf Lemma 2.1.} {\em A stochastic quadratic mapping $\:V\:$ is Bernstein 
iff the baric algebra $\:(\CA_V,s)\:$ is Bernstein in the sense}
$$
(x^2)^2\:=\:s^2(x)x^2.
\eqno (2.3)
$$

This lemma appeared first in [7] being written in the form
$\:\tilde{V}^2x\:=\:s^2(x)\tilde{V}x,\:$ but in [11] we already wrote (2.3).
Later Holgate [6] and author [13] considered the Bernstein algebras
by itself. (The term {\em Bernstein algebra} was introduced in [13].) 
In [17] (Sections 3.3 and 3.4) a part of the Bernstein algebras theory 
(over the field $\:\R )\:$ is presented 
in a form adapted to the Bernstein problem. Below we partially reproduce
it with addition of some new facts we need here.

In an arbitrary Bernstein algebra $\:(\CA,\gs )\:$  we have 
$$
\gs(xy)\:=\:\gs(x)\gs(y);\ \ \ (x^2)^2\:=\:\gs^2(x)x^2
\eqno (2.4)
$$
by definition. The first of these identities shows that the
subspace $\:\CB\:=\:{\rm Ker}\gs\:=\:\{\:x\: :\:x\:\in\:\CA,\ \ 
\gs(x)=0\:\}\:$ is an ideal (so-called {\em barideal}) in $\:\CA.\:$ The 
second one yields a construction of the idempotents in $\:\CA :\:$ if 
$\:\gs(x)\:=\:1\:$ then $\:e\:=\:x^2\:$ is an idempotent and 
$\:\gs(e)\:=\:1,\:$ so $\:e\:\neq\:0.\:$ (Conversely, if $\:e\:=\:e^2\:$
and $\:e\:\neq\:0\:$ then $\:\gs(e)\:=\:1.)\:$  

Given an idempotent $\:e\neq0,\:$ the linear operator $\:L_ey\:=\:2ey\:$ 
is a projection in $\:\CB\:$ hence, $\:\CB\:=\:U\:\oplus\:W\:$ where
$\:U\:=\:{\rm Im}L_e,\ \ W\:=\:{\rm Ker}L_e.\:$ Respectively,
$$
\CA\:=\:E\:\oplus\:U\:\oplus\:W
\eqno (2.5)
$$
where $\:E\:=\:{\rm Lin}\{e\},\:$ the linear span of $\:\{e\}.\:$ The
subspaces $\:E,\:U,\:W\:$ depend on $\:e\:$ but the dimensions 
$\:m-1\ {\rm and}\ \gd\:$ are invariant. 
The pair $\:(m,\gd)\:$ is called the {\em  type} of $\:\CA$.
Moreover $m={\rm rk}\CA$ is called the {\em rank} of
 $\:\CA$, and $\gd\:=\:{\rm def}A\: $
is called the {\em defect} of $\:\CA$.  Obviously,
$\:m\:+\:\gd\:=\:n\:=\:{\rm dim}\CA.\:$

The algebraic structure is reflected in (2.5) by the system of
inclusions:
$$
U^2\:\subset\:W,\ \ UW\:\subset\:U,\ \ W^2\:\subset\:U.
\eqno (2.6)
$$
Moreover, there is a series of identities connecting the variables 
$\:u\:\in\:U\:$ and $\:w\:\in\:W\:$ but we do not need this here.

If according to (2.5)
$$
x\:=\:\gs e\:\oplus\:u\:\oplus\:w
\eqno (2.7)
$$
then $\:\gs\:=\:\gs(x)\:$ and the corresponding decomposition of 
$\:x^2\:$ is
$$
x^2\:=\:\gs^2e\:\oplus\:(\gs u\:+\:2uw\:+\:w^2)\:\oplus\:u^2
\eqno (2.8)
$$
because of (2.6) and $\:2eu\:=\:L_eu\:=\:u,\ \ 2ew\:=\:L_ew\:=\:0.\:$

The simplest Bernstein algebras are the {\em constant algebras (c.a.)},
$\:x^2\:=\:c\gs^2(x)\:$ with $\:c\:\in\:\CA.\:$ All these algebras are of 
type $\:(1,n-1)\:$ and conversely, {\em every Bernstein algebra of type 
$\:(1,n-1)\:$ is constant}. An opposite simple example is the {\em unit
algebra (u.a.),} $\:x^2\:=\:\gs(x)x\:$ which is the only Bernstein algebra
of type $\:(n,0).\:$ {\em With $\:n\:\leq\:2\:$ every Bernstein algebra is
a c.a. or u.a.}

The evolutionary algebra $\:\CA_V\:$ is a c.a. (or u.a.) iff $\:V\:$ is  
constant (or identity) mapping.

The multiplication table for a c.a. is
$$
e_ie_k\:=\:c\ \ \ (1\leq i,k\leq n)
\eqno (2.9)
$$
and $\:xy\:=\:c\gs(x)\gs(y)\:$ for all $\:x,y.\:$

For the u.a. we have
$$
e_ie_k\:=\:\frac{e_i\:+\:e_k}{2}\ \ \ (1\leq i,k\leq n)
\eqno (2.10)
$$
and
$$
xy\:=\:\frac{\gs(y)x\:+\:\gs(x)y}{2}
\eqno (2.11)
$$
for all $\:x,y.\:$

The evolutionary algebra corresponding to the Hardy-Weinberg mapping is 
called the 
{\em Mendel algebra (M.a.)} This is a Bernstein algebra of type (2,1) 
with the multiplication table
$$
\left\{\begin{array}{lc}
e_1^2\:=\:e_1,\ e_2^2\:=\:e_2,\ e_1e_2\:=\:e_3,\ \ \ \ \ \ \ \ \ \ \ \ \ \ 
\ \ \ \ \ \ \ \ \ \ \ \ \ \ \ \ \ \ \ \ \ \ \ \ \ \ \ \ \ \ \ \  & 
\ \ \ \ \ \ \ \ \ \ \ \ \ \ \ \  (2.12)\\
\                                                               &        \\
e_1e_3\:=\:\frac{1}{2}(e_1+e_3),\ e_2e_3\:=\:\frac{1}{2}(e_2+e_3),\ \ \ 
\ \\ \ \ \ \ \ \ \ \ \ \ \ \ \ \ \ \ \ \ \ \ \ \ \ \ \ \ \ \ \ \ \ \ \ \ \ 
\ & \ \ \ \ \ \ \ \ \ \ \ \ \ \ \   (2.13)\\
\                                                                &       \\
e_3^2\:=\:\frac{1}{4}e_1\:+\:\frac{1}{4}e_2\:+\:\frac{1}{2}e_3.\ \ \ \ \ \ 
\ \ \ \ \ \ \ \ \ \ \ \ \ \ \ \ \ \ \ \ \ \ \ \ \ \ \ \ \ \ \ \ \ & 
\ \ \ \ \ \ \ \ \ \ \ \ \ \ \ \   (2.14) 
\end{array}
\right.
$$

Actually, (2.13) and (2.14) follow from (2.12) by the Bernstein property.

{\bf Proposition 2.2.} ([17], p.104) {\em Let $\:z_1\:$ and $\:z_2\:$ be 
nonzero idempotents in a Bernstein algebra 
$\:\CA\:$ and let $\:z_3=z_1z_2.\:$ Then}
$$
z_1z_3\:=\:\frac{z_1\:+\:z_3}{2},\ \ \ z_2z_3\:=\:\frac{z_2\:+\:z_3}{2}
\eqno (2.15)
$$
{\em and}
$$
z_3^2\:=\:\frac{1}{4}z_1\:+\:\frac{1}{4}z_2\:+\:\frac{1}{2}z_3.
\eqno (2.16)
$$

Therefore $\:Z\:=\:{\rm Lin}\{\:z_1,z_2,z_3\:\}\:$ is a Bernstein 
subalgebra which is isomorphic to the M.a. if dim$Z\:=\:3,\:$ i.e. if
$\:z_1,z_2,z_3\:$ are linearly independent. If dim$Z\:=\:2\:$ (i.e.
$\:z_1\neq z_2\:$ and $\:z_3\:\in\:{\rm Lin}\{\:z_1,z_2\:\})\:$ then 
$\:Z\:$ is the u.a.

It is very useful for our purposes to introduce a new commutative 
multiplication,
$$
R(x,y)\:=\:2xy\:-\:\gs(y)x\:-\:\gs(x)y
\eqno (2.17)
$$
in a Bernstein algebra $\:\CA.\:$ Letting $\:x\circ y\:$ for the unit 
multiplication (2.11) we obtain
$$
R(x,y)\:=\:2(xy\:-\:x\circ y)
\eqno (2.18)
$$
So $\:R\:$ measures a deviation of the given algebra from the u.a.,
$$
R(x,y)\:=\:0\ \Leftrightarrow\ xy\:=\:x\circ y
\eqno (2.19)
$$

Obviously, $\:R(x,x)\:=\:0\:$ iff $\:x^2\:=\:\gs(x)x,\:$ in 
particular, $\:R(x,x)\:=\:0\:$ for all idempotents $\:x.\:$

Note that the subalgebras are the same for $\:R(x,y)\:$ and $\:xy\:$ 
(including the non-Bernstein ones in $\:(\CA,\gs),\:$ i.e. the 
subalgebras of the barideal $\:\CB ).\:$

{\bf Lemma 2.3.} {\em Any four idempotents $\:z_1,\:z_2.\:z_3,\:z_4\:$
in a Bernstein algebra satisfy the relation}
$$
R(z_1z_2,z_3z_4)\:+\:R(z_1z_3,z_2z_4)\:+\:R(z_1z_4,z_2z_3)\:=\:0.
\eqno (2.20)
$$

{\bf Proof.} It is trivial if one of $\:z_i\:$ is zero. If all of   
them are nonzero we insert
$$
x\:=\:\sum\limits_{i=1}^4\xi_iz_i,\ \ \ \gs(x)\:=\:\sum\limits_{i=1}^4\xi_i
$$
into the identity $\:(x^2)^2\:=\:\gs^2(x)x^2\:$ and then compare the
coefficients at the monomial $\:\xi_1\xi_2\xi_3\xi_4.\:$ This immediately
leads to (2.20). 
$\ \ \ \ \ \ \ \ \ \ \ \ \ \ \ \ \ \ \ \ \ \ \ \ \ \ \ \ \ \ \ \ \ \ \ \ 
\ \ \ \ \ \ \ \ \ \ \ \ \ \ \ \ \ \ \ \ \ \ \ \ \Box $

{\bf Corolary 2.4.} {\em For any three idempotents $\:z_1,z_2,z_3\:$}
$$
2R(z_1z_2,z_1z_3)\:+\:R(z_1,z_2z_3)\:=\:0
\eqno (2.21)
$$

{\bf Proof.} Take $\:z_4\:=\:z_1\:$ in (2.20).
$\ \ \ \ \ \ \ \ \ \ \ \ \ \ \ \ \ \ \ \ \ \ \ \ \ \ \ \ \ \ \ \ \ \ \ \ 
\ \ \ \ \ \ \ \  \ \ \ \ \ \ \ \ \ \ \ \ \ \ \ \ \ \ \Box\ $

We also come back to (2.15) and (2.16) setting $\:z_3\:=\:z_1\:$ or 
$\:z_3\:=\:z_2\:$ in (2.21).

For some further constructions we need 

{\bf Corollary 2.5.} {\em If $\:z_1,z_2,w_1,w_2\:$ are nonzero 
idempotents such that}
$$
R(z_i,w_j)\:=\:0\ \ \ (i,j\:=\:1,2)
\eqno (2.22)
$$
{\em then}
$$
R(z_1z_2,w_1w_2)\:=\:-\frac{R(z_1,z_2)\:+\:R(w_1,w_2)}{2}.
\eqno (2.23)
$$

{\bf Proof.} It follows from (2.20) and (2.19) (by assumption (2.22)) that
$$
R(z_1z_2,w_1w_2)\:=\:-R(z_1w_1,z_2w_2)\:-\:R(z_1w_2,z_2w_1)\:=
$$
$$
-R(z_1\circ w_1,z_2\circ w_2)\:-\:R(z_1\circ w_2,z_2\circ w_1)\:=
$$
$$
=-\frac{1}{4}[R(z_1\:+\:w_1,\:z_2\:+\:w_2)\:+\:
R(z_1\:+\:w_2,\:z_2\:+\:w_1)]
$$
which can be reduced to (2.23) by (2.22).
$\ \ \ \ \ \ \ \ \ \ \ \ \ \ \ \ \ \ \ \ \ \ \ \ \ \ \ \ \ \ \ \ 
\ \ \ \ 
\ \ \ \ \ \ \ \ \ \ \ \ \ \ \ \ \ \ \ \ \ \ \ \ \ \ \ \ \ \ \ \ \ \ \Box\ $

{\bf Corollary 2.6.} {\em If $\:z_1,z_2,w\:$ are nonzero idempotents such
that}
$$
R(z_i,w)\:=\:0\ \ \ (i\:=\:1,2)
\eqno (2.24)
$$
{\em then}
$$
R(z_1z_2,w)\:=\:-\frac{1}{2}R(z_1,z_2).
\eqno (2.25)
$$

Therefore $\:R(z_1z_2,w)\:$ is independent of $\:w.\:$

{\bf Proof.} Take $\:w_1\:=\:w_2\:=\:w\:$ in (2.23). 
$\ \ \ \ \ \ \ \ \ \ \ \ \ \ \ \ \ \ \ \ \ \ \ \ \ \ \ \ \ \ \ \ \ \ \ \ 
\ \ \ \ \ \ \ \ \ \ \ \ \ \ \ \Box\ $

In our context the most important baric algebras are {\em regular} ones.
By one of many equivalent definitions, the {\em regularity} 
of a baric algebra $\:(\CA ,\gs )\:$ means that
$\:xy\:$ only depends on values $\:f(x)\:$ and $\:f(y)\:$ where $\:f\:$ 
runs over all {\em invariant} linear forms. The {\em invariance} of 
$\:f\:$ means that $$
\gs(x)\:=\:1\ \Rightarrow\ f(x^2)\:=\:f(x),
\eqno (2.26)
$$
or equivalently, 
$$
f(xy)\:=\:\frac{\gs(y)f(x)\:+\:\gs(x)f(y)}{2}\
\eqno (2.27)
$$
which in turn can be written as
$$
f(R(x,y))\:=\:0.
\eqno (2.28)
$$
For the evolutionary algebras $\:\CA_V\:$ the invariance of $\:f\:$ is the 
same as for $\:V,\:$
i.e. $\:f(Vx)\:=\:f(x).\:$ Therefore, for any $\:(\CA ,\gs )\:$ 
we can use the notation $\:J\:$
for the space of all linear invariant forms. Obviosly,
$\:\gs\:\in\:J,\:$ so dim$J\:\geq\:1.\:$

{\em An evolutionary operator $\:V\:$ is regular iff
the algebra $\:\CA_V\:$ is regular.} 

Note that the invariant faces of
$\:\gD^{n-1}\:$ are just such that their linear spans are subalgebras in
$\:\CA_V,\:$ i.e. they are coordinate subalgebras. Let us say that 
$\:\CA_V\:$ is {\em normal} if $\:V\:$ is so (see [17], Section 3.9 for a
more algebraic treat of this notion). Respectively, $\:\CA_V\:$ is said
to be {\em ultranormal} if $\:V\:$ is so, i.e. all coordinate subalgebras
are normal. Thus, we are going to prove

{\bf The Main Theorem}.  
{\em Every ultranormal stochactic Bernstein algebra is regular.}

For this goal we need some 
regularity criteria for the Bernstein algebras. 

Certainly, any 
regular algebra is Bernstein. This easily follows from definitions or 
from the identity
$$
x^2y\:=\:\gs(x)xy
\eqno (2.29)
$$
characterizing the regularity ([13]; [17], 
Theorem 3.3.6).(By the way, such a characterization shows that all 
subalgebras of a regular algebra are regular).

Note that dim$J\leq m\:$ for any Bernstein of rank$\:m.\:$
 
{\bf Theorem 2.7.} {\em For any Bernstein algebra of rank$\:m\:$ the
following conditions are equivalent:}

1) {\em the algebra is regular;}

2) dim$J\:=\:m;$

3) $UW\:+\:W^2\:=\:0;$

4) $UW\:=\:0\:$ {\em and $\:W^2\:=\:0,\:$ so that (2.8) takes the form}
$$
x^2\:=\:\gs^2e\:\oplus\:\gs u\:\oplus\:u^2.
\eqno (2.30)
$$
(see [7]; [17], Theorems 3.3.4, 3.4.15 and 3.4.17).

A very important consequence of this criterion is

{\bf Corollary 2.8.} {\em A Bernstein algebra is regular if (2.29) holds
for all of $\:x\:$ and for a fixed idempotent $\:y\:=\:e\:\neq\:0.\:$}

{\bf Proof.} Insert $\:x\:$ and $\:x^2\:$ from (2.7) and (2.8) into (2.29). 
We get 
$$
\gs^2e\:\oplus\:\frac{1}{2}(\gs u\:+\:2uw\:+\:w^2)\:=\:\gs^2e\:\oplus\:
\frac{1}{2}\gs u
$$
because of $\:e^2\:=\:e,\ \ eu\:=\:\frac{1}{2}u,\ \ ew\:=\:0.\:$ Hence
$\:uw\:+\:\frac{1}{2}w^2\:=\:0.\:$ The algebra is regular by part 3) of
Theorem 2.7. 
$\ \ \ \ \ \ \ \ \ \ \ \ \ \ \ \ \ \ \ \ \ \ \ \ \ \ \ \ \ \
\ \ \ \ \ \ \ \ \ \ \ \ \ \ \ \ \ \ \ \ \ \ \ \ \
\ \ \ \ \ \ \ \ \ \ \ \ \ \ \ \ \ \ \ \ \ \ \ \ \ \ \ \ \Box\ $ 

It is conveniente to formulate Corollary 2.8 in a coordinate from.

{\bf Corollary 2.9.} {\em Let a Bernstein algebra $\:\CA\:$ is a 
linear span of 
a system of vectors $\:\{\:v_i\:\}_1^l,\ \ \gs(v_i)\:=\:1\ 
\ (1\:\leq\:i\:\leq\:l).\:$ If there exists an idempotent$\:e\neq 0\:$ 
such that } $$
(v_iv_k)e\:=\:\frac{v_ie\:+\:v_ke}{2}\ \ (1\:\leq\:i,k\:\leq\:l)
\eqno (2.31)
$$
{\em then the algebra is regular.}

{\bf Proof.} Any vector $\:x\:\in\:\CA\:$ is  
$\:x\:=\:\sum\xi_iv_i.\:$ Then (2.31) implies $\:x^2e\:=\:\gs(x)xe\:$ since
$\:\gs(x)\:=\:\sum\xi_i.\:$ 
$\ \ \ \ \ \ \ \ \ \ \ \ \ \ \ \ \ \ \ \ \ \ \ \ \ \ \ \ \ \ \ \ \ \ \ \ \ 
\ \ \ \ \ \ \ \ \ \ \ \ \ \ \ \ \ \ \ \ \ \ \ 
\ \ \ \ \ \ \ \ \ \ \ \ \ \ \ \ \ \ \ \ \ \ \ \ \ \ \ \ \ \ \ \Box\ $

Note that (2.31) can be rewriten as
$$
R(v_i,v_k)e\:=\:0\ \ \ (1\:\leq\:i,k\:\leq\:l).
\eqno (2.32)
$$
Similarly, (2.29) can be rewritten as $\:R(x,x)y\:=\:0.\:$ This identity 
is equivalent to $\:R(x,z)y\:=\:0\:$ which is formally a more general one.

The concrete examples of regular algebras are c.a.,  u.a.,  M.a.
By Proposition 2.2 any pair of idempotents $\:z_1,z_2\ \ (z_1\:\neq\:z_2)\:$
generates either the (2-dimensional) u.a. or the M.a. Using more idempotents 
one can inductively construct some other regular subalgebras in a Bernstein
algebra. 

{\bf Proposition 2.10.} {\em Let $\:\{\:z_i\:\}_1^{\nu}\:$ be a family of 
idempotents such that the subspace $\:L\:=\:{\rm Lin}\:\{\:z_iz_k\:\}_
{i,k=1}^{\nu}\:$ is a regular subalgebra. Then}

1) {\em if $\:w\:$ is an idempotent such that}
$$
R(z_j,w)\:=\:0\ \ \ (1\:\leq\:j\:\leq\:\nu)
\eqno (2.33)
$$
{\em then $\:L[w]\:=\:L\:\oplus\:{\rm Lin}\:\{\:w\:\}\:$ is a regular
subalgebra;}

2) {\em if $\:w_1\:$ and $\:w_2\:$ are idempotents such that}
$$
R(z_i,w_j)\:=\:0\ \ \ (1\leq i\leq\nu ;\ j=1,2)
\eqno (2.34)
$$
{\em and}
$$
R(w_1,w_2)z_1\:=\:0
\eqno (2.35)
$$
{\em then $\:L[w_1,w_2]\:=\:L\:\oplus\:{\rm 
Lin}\:\{\:w_1,w_2,w_1w_2\:\}\:$ is a regular subalgebra.}

{\bf Proof.} 1) By Corollary 2.6 all $\:R(z_iz_k,w)\:\in\:L,\:$ so
$\:(z_iz_k)w\:\in\:L[w]\:$ 
hence, $\:L[w]\:$ is a subalgebra. This is the
linear span of $\:w\:$ and all of $\:z_iz_k\ \ (1\leq i,k\leq\nu ),\ z_1\:$
among them. By Corollary 2.9 with $\:e=z_1\:$ 
$\:L[w]\:$ is regular. Indeed, by Corollary 2.6
$$
z_1R(z_iz_k,w)\:=\:-\frac{1}{2}z_1R(z_i,z_k)\ \ (1\leq i,k\leq\nu )
$$
which equals zero because $\:L\:$ is regular.

2) We already know that $\:L[w_1],\:L[w_2]\:$ and 
Lin$\{\:w_1,w_2,w_1w_2\:\}\:$ are regular subalgebras. Furthermore, all
$\:R(z_iz_k,w_1w_2)\:\in\:L[w_1,w_2]\:$ by Corollary 2.5. So,
$\:(z_iz_k)(w_1w_2)\:\in\:L[w_1,w_2].\:$ We see that $\:L[w_1,w_2]\:$ is
a subalgebra. It is regular by Corollary 2.5 and the regularity of $\:L\:$
imply
$$
z_1R(z_iz_k,w_1w_2)\:=\:-\frac{1}{2}z_1R(w_1,w_2)\:=\:0\ \ 
(1\leq i,k\leq\nu ),
$$
and, moreover,
$$
R(w_j,w_1w_2)\:=\:0\ \ (j=1,2),\ \ 
R(w_1w_2,w_1w_2)\:=\:-\frac{1}{2}R(w_1,w_2)
$$
by Proposition 2.2.
$\ \ \ \ \ \ \ \ \ \ \ \ \ \ \ \ \ \ \ \ \ \ \ \ \ \ \ \ \ \ \ \ \ \ \ \ 
\ \ \ \ \ \ \ \ \ \ \ \ \ \ \ \ \ \ \ \ \ \ \ \ \ \ \ \ \ \ \ \ \ \ \ \ \
\ \ \ \ \ \ \Box\ $

In conclusion we formulate a sufficient regularity condition in terms of 
type $\:(m,\gd)\:$ ([10]; [17], Corollaries 3.4.28 and 3.4.29).

{\bf Theorem 2.11.} {\em Let a Bernstein algebra $\:\CA\:$ of type 
$\:(m,\gd)\:$ be nuclear, i.e. $\:\CA^2\:=\:\CA.\:$ Then it is regular if
$\:m\:\leq\:3\:$ or $\:\gd\:\leq\:1,\:$ or}
$$
\gd\:\geq\:\frac{(m-1)(m-2)}{2}\:+\:1.
\eqno (2.36)
$$

{\bf Corollary 2.12.} {\em Every nuclear Bernstein algebra of dimension 
$\:n\leq 5\:$ is regular.}

There exists a nonregular nuclear Bernstein algebra of type (4,2) (see [17],
p.102). The question whether {\em every stochastic nuclear
Bernstein algebra is regular} (the Conjecture 5.7.16 from [17]) is still 
open. However, we have

{\bf Theorem 2.13.} {\em With $\:m\:\leq\:2\:$ or $\:\gd\:\leq\:1\:$
(in particular, with $\:n\:\leq\:4\:$ ) every normal stochastic Bernstein 
algebra is regular.}

Finally, we have such a part of the Main Theorem.

{\bf Theorem 2.14} ( [16]; [17], Section 5.7). {\em With $\:m\:\leq\:3\:$ or
$\:\gd\:\leq\:1\:$ (with $\:n\:\leq\:5,\:$ in particular) every ultranormal
stochastic Bernstein algebra $\:\CA\:$ is regular.}

Its proof is based on a combinatorial topology structure which is  
induced by invariant faces on the set Im$V\:$ where $\:V\:$ is the
corresponding quadratic mapping, $\:\CA_V\:=\:\CA .\:$ 
We develop this approach 
in Section 4. 

\section {Normal stochastic regular algebras} 

$\ \ \ \ \ \ \ \ $ The complete solution of the Bernstein problem for the 
normal regular algebras is given by the following theorem which is an 
algebraic reformulation of Theorems 4.3.9 and 4.6.1 from [17].

{\bf Theorem 3.1.} {\em Every normal stochastic regular algebra $\:\CA\:$ 
of type $\:(m,\gd)\:$ is one of two following ones.}

1) {\em Up to enumeration of the canonical basis $\:\{\:e_i\:\}^n_1\:$ 
the vectors $\:e_1,\ ...\:,e_m\:$ are idempotents. Their products are}
$$
e_{i_j}e_{k_j}\:=\:\ga_je_{i_j}\:+\:\gb_je_{k_j}\:+\:\gamma_je_{m+j}
\eqno (3.1)
$$
{\em for some distinct pairs $\:(i_j,k_j)\:$ with $\:1\:\leq\:i_j\:<\:k_j\:
\leq\:m,\ \ 1\:\leq\:j\:\leq\:\gd,\:$ and}
$$
e_ie_k\:=\:\frac{e_i\:+\:e_k}{2}
\eqno (3.2)
$$
{\em for all remaining pairs} $\:(i,k)\ \ (1\leq i<k\leq m,\ \ i\neq i_j\:$ 
{\em or} $\: k\neq k_j).\:$ {\em In (3.1)} $\:\ga_j\:\geq\:0,\ 
\gb_j\:\geq\:0,\ \ \ga_j\:+\:\gb_j\:+\:\gamma_j\:=\:1.$ 

{\em Furthermore,}
{\em for $\:1\:\leq\:j\:\leq\:\gd\:$ and $\:1\:\leq\:i\:\leq\:m\:$}
$$
e_ie_{m+j}\:=\:c_je_ie_{i_j}\:+\:\overline{c_j}e_ie{k_j}
\eqno (3.3)
$$
{\em where $\:0\:<\:c_j\:<\:1;\ \ \overline{c_j}\:=\:1-c_j\:$ and }
$$
\ga_j\:+\:c_j\gamma_j\:=\:\gb_j\:+\:\overline{c_j}\gamma_j\:=\:\frac{1}{2}.
\eqno (3.4)
$$
Finally,
$$
e_{m+j}e_{m+l}\:=\:c_jc_le_{i_j}e_{i_l}\:+\:c_j\overline{c_l}e_{i_j}e_{k_l}
\:+\:\overline{c_j}c_le_{k_j}e_{i_l}\:+\:\overline{c_j}\:\overline{c_l}
e_{k_j}e_{k_l}
\eqno (3.5)
$$
{\em for $\:1\:\leq\:j,l\:\leq\:\gd.\:$ }

2){\em All basis vectors $\:\{\:e_i\:\}^n_1\:$ are idempotents. The number
$\:n\:=\:{\rm dim}\CA\:$ is composite, $\:n\:=\:\nu\overline{\nu}\:$ with
$\:\nu\:\geq\:2,\ \overline{\nu}\:\geq\:2.\:$}
{\em The basis can be enumerated as $\:\{\:e_{ik}\: :\:1\leq i\leq\nu,\ 
1\leq k\leq\overline{\nu}\:\}\:$ in a way such that }
$$
e_{gh}e_{ik}\:=\:\frac{e_{gk}\:+\:e_{hi}}{2}
\eqno (3.6)
$$
{\em for all pairs $\:(g,h)\:$ and $\:(i,k).\:$}

The first case is the {\em elementary gene structure (e.g.s.),}
the second one is the {\em nonelementary gene structure (n.e.g.s.).}
Both of them allow a natural genetical interpretation (see [17], p.p. 189,
207). E.g.s. is {\em continual,} i.e. multiparametric.
The independent parameters are $\:\ga_j,\gb_j\ \ (1\leq j\leq\gd)\:$ so
that the manifold of all these algebras is
$\:2\gd$-dimensional.
N.e.g.s. is {\em discrete,} i.e. 0-dimensional.

The unit algebra (u.a.) has e.g.s. In this case $\:\gd\:=\:0\:$ and there 
is no $\:\:e_{m+j},\:$ no pairs $\:(\:i_j,k_j\:),\:$ so that (3.2) is the 
complete multiplication table. It is the only case with e.g.s. when all 
$\:e_i\ \ (1\leq i\leq n)\:$ are idempotents.

The constant algebra (c.a.) is normal in the only case $\:n=1\:$ but then 
it is also u.a.

The simplest nontrivial situation is 3-dimensional.

{\bf Example 3.2.} For $\:m=2\:$ and $\:\gd=1\ \ (n=3)\:$ we have e.g.s.
$$
e_1^2\:=\:e_1,\ \ e_2^2\:=\:e_2,\ \ 
e_1e_2\:=\:\ga e_1\:+\:\gb e_2\:+\:\gamma e_3
\eqno (3.7)
$$
where $\:\ga\geq0,\ \gb\geq0,\ \gamma >0,\ \ga+\gb+\gamma\:=\:1.\:$ 
Furthermore,
$$
e_3e_1\:=\:ce_1^2\:+\:\overline{c}e_1e_2,
$$
so that
$$
e_3e_1\:=\:(c+\overline{c}\ga)e_1\:+\:\overline{c}\gb e_2\:+\:\overline{c}
\gamma e_3
\eqno (3.8)
$$
and similarly,
$$
e_3e_2\:=\:ce_1e_2\:+\:\overline{c}e_2^2\:=\:c\ga e_1\:+\:(c\gb 
+\overline{c})e_2\:+\:c\gamma e_3
\eqno (3.9)
$$
with $\:0\:<\:c\:<\:1,\ \ \overline{c}\:=\:1\:-\:c\:$ and
$$
\ga\:+\:\gamma c\:=\:\gb\:+\:\gamma\overline{c}\:=\:\frac{1}{2}.
\eqno (3.10)
$$
Finally,
$$
e_3^2\:=\:c^2e_1^2\:+\:2c\overline{c}e_1e_2\:+\:\overline{c}^2e_2^2,
$$
so that
$$
e_3^2\:=\:(c^2\:+\:2c\overline{c}\ga)e_1\:+\:(2c\overline{c}\gb\:+\:
\overline{c}^2)e_2\:+\:2c\overline{c}\gamma e_3.
\eqno (3.11)
$$

The Mendel algebra (M.a.) is just the case $\:\ga\:=\:\gb\:=\:0,\ \ 
\gamma\:=\:1,\ \ c\:=\:\frac{1}{2},\:$ so $\:\overline{c}\:=\:\frac{1}{2}.\:$
M.a. is the point $\:M(0,0)\:$ of the 2-dimensional manifold of algebras
$\:M(\ga,\gb)\:$ given by (3.7)-(3.11). Therefore $\:M(\ga,\gb)\:$ can be
called an {\em extended Mendel algebra (e.M.a.).}

Let us emphasize that all the algebras $\:M(\ga,\gb)\:$ are isomorphic to the
M.a. $\:M(0,0)\:$ since (3.7) with $\:\gamma >0\:$ allows us to change the 
basis $\:\{\:e_i\:\}_1^3\:$ for $\:\{\:e_1,e_2,e_1e_2\:\}\:$ and then the 
algebra turns into the M.a. Proposition 2.2. However, we must distinguish
the algebras $\:M(\ga,\gb)\:$ at the fixed basis $\:\{\:e_i\:\}_1^3.\:$ 
Biologically, it is necessary because $\:e_i\:$ are the types themselves 
but their linear combinations have no such interpretation; respectively, the
coefficients $\:\ga,\gb,\gamma,\:$ etc. are the probabilities of types in the
offsprings generation.

The term {\em e.M.a.} can be also used for any algebra of type $\:(m,\gd)\:$
given by (3.1)-(3.5). The corresponding evolutionary operator is an {\em 
extended Hardy-Weibnerg mapping.} Actually, it is 
$$
x_i^{'}\:=\:p_i^2\:+\:2\sum\limits_{k\neq i}\gt_{ik}p_ip_k\ \ \ 
(1\:\leq\:i\:\leq\:m)
\eqno (3.12)
$$
and
$$
x_{m+j}^{'}\:=\:2\gamma_jp_{i_j}p_{k_j}\ \ \ (1\:\leq\:j\:\leq\:\gd)
\eqno (3.13)
$$
where 
$$
p_i\:=\:x_i\:+\:\sum\limits_{j=1}^{\gd}\pi_{ij}x_{m+j},
\eqno (3.14)
$$
$\:\gt_{i_jk_j}\:=\:\ga_j,\ \ \gt_{k_ji_j}\:=\:\gb_j\:$ and all 
remaining $\:\gt_{ik}\:=\:\frac{1}{2};\ \ \pi_{i_jj}\:=\:c_j;\ \ 
\pi_{k_jj}\:=
\:\overline{c_j}\:$ and all remaining $\:\pi_{ij}\:=\:0\:$ (cf. [17],
Theorem 4.3.9).

The set $\:G\:=\:\{\:p_i\:\}_1^m\:$ is just the canonical basis of the cone
$\:C\:$ of all nonnegative invariant linear forms. Obviously, this set is 
linearly independent and
$$
\sum\limits_{i=1}^mp_i\:=\:s
\eqno (3.15)
$$
The cone $\:C\:$ is minihedral in the case of e.g.s.

{\bf Example 3.3.} For $\:\nu\:=\:\overline{\nu}\:=\:2\ \ (n=4)\:$ we
have the symplest n.e.g.s. which is actually the {\em quadrille algebra 
(q.a.)} corresponding to the quadrille mapping (1.9)-(1.10). A more 
natural labeling in this case is $\:x_1\:\equiv\:x_{11},\ \ x_2\:\equiv\:
x_{22},\ \ x_3\:\equiv\:x_{12}\:$ and $\:x_4\:\equiv\:x_{21}.\:$ As a result
$$
x_{ik}^{'}\:=\:p_iq_k
\eqno (3.16)
$$
where $\:p_i\:$ are the sums over rows of the matrix 
$\:X\:\equiv\:(x_{ik})\:$ and $\:q_k\:$ are the sums over columns,
$$
p_i\:=\:\sum\limits_kx_{ik},\ \ \ q_k\:=\:\sum\limits_ix_{ik}.
\eqno (3.17)
$$
In a matrix form (3.6) is 
$$
X^{'}\:=\:p\:\otimes\:q
\eqno (3.18)
$$
where $\:p\:$ is the column $\:(p_i)\:$ and $\:q\:$ is the row $\:(q_k).\:$

The same formulae (3.16)-(3.18) take place in general, i.e. for any $\:n\:=\:
\nu\overline{\nu}\:$ with $\:X\:=\:(x_{ik}\: :\:1\leq i\leq\nu,\ \ 1\leq 
k\leq\overline{\nu}).\:$ Any such a mapping is called an {\em extended 
quadrill mapping} and the corresponding algebra is an 
{\em extended quadrille algebra (e.q.a.)}.

The set $\:G\:=\:\{\:p_i\:\}_1^{\nu}\:\bigcup\:\{\:q_k\:\}_1^{\overline{\nu}}
\:$ is the canonical basis of the cone $\:C\:$ in this case. Now this set 
is linearly dependent since
$$
\sum\limits_{i=1}^{\nu}p_i\:=\:\sum\limits_{k=1}^{\overline{\nu}}q_k\ \ 
\ (=s). 
\eqno (3.19)
$$
Thus, the cone $\:C\:$ is not minihedral in the case of n.e.g.s.

Besides (3.19), there is no linear dependence in 
$\:G.\:$ Therefore the type $\:(m,\gd)\:$of the e.q.a. is
$$
m\:=\:\nu\:+\:\overline{\nu}\:-\:1,\ \ \ \gd\:=\:(\nu -1)(\overline{\nu} -1).
\eqno (3.20)
$$
In particular, the q.a. is of type (3,1).

{\bf Corollary 3.4.} {\em Every coordinate subalgebra of e.M.a. is also 
e.M.a. Every coordinate subalgebra of e.q.a is also e.q.a. or unit.}

{\bf Proof.} If follows from (3.5) with $\:l=j\:$ that every coordinate 
subalgebra of an e.M.a. containing $\:e_{m+j}\ \ (1\leq j\leq\gd)\:$ must
contain both of the idempotents $\:e_{i_j}\:$ and $\:e_{k_j}.\:$ The converse
is also true by (3.1). Thus, any coordinate subalgebra of the e.M.a. is
the linear span of the union of the subset $\:F\:\subset\:\{\:e_i\:\}_1^m\:$ 
with all of
$\:\{\:e_{m+j}\: :\:e_{i_j},e_{k_j}\:\in\:F\:\} .\:$ Obviously, it is an 
e.M.a. as well.

The case of e.q.a. is similar (even simpler).
$\ \ \ \ \ \ \ \ \ \ \ \ \ \ \ \ \ \ \ \ \ \ \  
\ \ \ \ \ \ \ \ \ \ \ \ \ \ \ \ \ \ \ \ \ \ \ \ \ \Box\ $
 
Note that {\em any e.M.a. is regular and normal} because of (3.12)-(3.14) 
where the pairs $\:(i_j,k_j)\:$ are distinct and the restrictions 
$\:0<c_j<1,\ \ \gamma_j>0\:$ are fulfilled.
{\em Any e.q.a. is also regular and normal} because of (3.16)-(3.17).
By Theorem 3.1 and Corollary 3.4 we get

{\bf Corollary 3.5.} {\em Every normal stochastic regular algebra is 
ultranormal.}

In addition, we have

{\bf Corollary 3.6.} {\em Every normal stochastic regular algebra is 
nuclear.}

{\bf Proof.} In the case of e.M.a.
$$
{\rm Lin}\{\:e_ie_k\:\}_{i,k=1}^n\:=\:
{\rm Lin}\{\:e_i\:\}_1^m\:\bigcup\:{\rm 
Lin}\{\:e_{i_j}e_{k_j}\:\}_1^{\gd}\:=\:\CA
$$
because the second set in the union can be changed for  
Lin$\{\:e_{m+j}\:\}_1^{\gd}\:$ using (3.1) with $\:\gamma_j>0.\:$ 
In the case of e.q.a. $\:\CA^2\:=\:\CA\:$
because of $\:e_i^2\:=\:e_i\:$ for all $\:i,\ 1\leq i\leq n.\:$
$\ \ \ \ \ \ \ \ \ \ \ \ \ \ \ \ \ \ \ \ \ \ \ \ \ \ \ \Box\ $

\section{Offspring subalgebras}
$\ \ \ \ \ \ \ \ $
Recall that for any vector $\:x\:\in\:\R^n,\:$
$$
x\:=\:\sum\limits_{j=1}^nx_ie_i,
$$
its {\em support} is defined as
$$
{\rm supp}x\:=\:\{\:e_i\: :\:x_i\neq 0\:\}
$$
so that
$$
x\:=\:\sum\limits\:\{\:x_ie_i\: :\:e_i\:\in\:{\rm supp}x\:\} .
\eqno (4.1)
$$
Obviously, supp$x\:\neq\:\emptyset\:$ for $\:x\:\neq\:0\:$ and
$$
{\rm supp}(\gl x)\:=\:{\rm supp}x\ \ \ (\gl\neq 0).
\eqno (4.2)
$$
If $\:x\geq 0\:$ and $\:y\geq 0\:$ then
$$
{\rm supp}(x+y)\:=\:{\rm supp}x\:\bigcup\:{\rm supp}y.
\eqno (4.3)
$$

We say that an algebra $\:\CA\:$ with the underlying space $\:\R^n\:$ is
{\em nonnegative} if
$$
x\:\geq\:0\:\&\:y\:\geq\:0\:\Rightarrow\:xy\:\geq\:0
\eqno (4.4.)
$$
or, equivalently, its structure constants are nonnegative. Every 
stochastic algebra is so.

{\bf Lemma 4.1.} {\em In any nonnegative algebra $\:\CA\:$ for any 
$\:x\geq 0\:$}
$$
{\rm supp}(x^2)\:=\:\bigcup\:\{\:{\rm supp}(e_ie_k)\: :\:e_i,e_k\:\in\:
{\rm supp}x\:\}.
\eqno (4.5)
$$

Thus, supp$(x^2)\:$ only depends on supp$x\ \ (x\geq 0).\:$

{\bf Proof.} It follows from (4.1) that
$$
x^2\:=\:\sum\limits\:\{\:x_ix_ke_ie_k\: :\:e_i,e_k\:\in\:{\rm supp}x\:\}
$$
and then (4.5) follows from (4.2) and (4.3).
$\ \ \ \ \ \ \ \ \ \ \ \ \ \ \ \ \ \ \ \ \ \ \ \ \ \ \ \ \ \ \ \ 
\ \ \ \ \ \ \ \ \ \ \ \ \ \ \ \ \ \ \ \ \ \ \ \ \ \ \ \Box\ $

Note that 
$$
{\rm supp}(e_ie_k)\:=\:\{\:e_j\: :\:p_{ik,j}>0\:\}.
$$
Biologically, supp$(e_ie_k)\:$ is the set of all types (characters) 
really presented in offsprings whose parental types are $\:e_i\:$ and 
$\:e_k.\:$ If $\:e_i^2\:=\:e_i\:$ the type $\:e_i\:$ is {\em nonsplitting}
in the sense that all its offspring are of the same type.

For any family $\:F\:\subset\:\{\:e_i\:\}_1^n\:$ we define its {\em 
offspring set}
$$
F^{'}\:=\:\bigcup\:\{ {\rm supp}(e_ie_k):\ \ e_i,e_k\:\in\:F\:\} .
$$
Vice versa $\:F\:$ is the {\em parental set of} $\:F^{'}.\:$ In the most
important case $\:F\:$ consists of some idempotents. Then 
$\:F^{'}\:\supset\:F.\:$

{\bf Lemma 4.2.} Lin$F\:$ {\em is a subalgebra iff $\:F^{'}\:\subset\:F.\:$}

{\bf Proof.} $\:(\A e_i,e_k\:\in\:F\: :\:e_ie_k\:\in\:{\rm Lin}F)\:
\Leftrightarrow\:(\A e_i,e_k\:\in\:F\: :\:{\rm 
supp}(e_ie_k)\:\subset\:F)\:\Leftrightarrow\:(\:\{\:\bigcup\: {\rm supp}
(e_ie_k)\: :\:e_i,e_k\:\in\:F\:\}\:\subset\:F)\:\Leftrightarrow\:
(F^{'}\:\subset\:F)\:$ by Lemma 4.1.      

{\bf Corollary 4.3.} Lin$\:\{\:{\rm supp}x\:\}\:$ {\em is a subalgebra for 
any idempotent $\:x\:\geq\:0.\:$ }

{\em Henceforth we only consider a stochastic Bernstein algebra 
$\:(\CA,s),\:$ so $\:\CA\:=\:\CA_V\:$ where $\:V\: 
:\:\gD^{n-1}\:\rightarrow\:\gD^{n-1}\:$ is a stochastic quadratic mapping,
$\:V^2\:=\:V.\:$ }

{\bf Lemma 4.4.} {\em $\:F^{''}\:=\:F^{'}\:$ for any family $\:F\:\subset\:
\{\:e_i\:\}_1^n.\:$ }

In this sense there are no characters coming from the offspring to their  
offspring but not originating from their parents.

{\bf Proof.} Let
$$
x\:=\:\sum\limits\:\{\:e_i\: :\:e_i\:\in\:F\:\},
$$
so that supp$x\:=\:F.\:$ Then supp$(x^2)\:=\:F^{'}\:$ and
$$
F^{''}\:=\:{\rm supp}(x^2)^2\:=\:{\rm supp}[s^2(x)x^2]\:=\:{\rm 
supp}(x^2)\:=\:F^{'}.
$$

{\bf Corollary 4.5.} Lin$F^{'}\:$ {\em is a subalgebra.}

We call this the {\em offspring subalgebra} of the parental set $\:F.\:$
This construction plays a very important role in sequel.

For example, $\:M (\ga,\gb)\:$ is the offspring subalgebra of the set
$\:F\:=\:\{\:e_1,e_2\:\}\:$ with $\:e_1^2\:=\:e_1,\ e_2^2\:=\:e_2.\:$ 
Moreover, e.M.a. end e.q.a. are both the offspring subalgebras of some 
families consisting of the basis idempotents. Such a family is 
$\:F\:=\:\{\:e_i\:\}_1^m\:$ for e.M.a. and 
$\:F\:=\:\{\:e_{g1}\:\}_{g=1}^{\nu}\ 
\bigcup\:\{\:e_{1h}\:\}_{h=1}^{\overline{\nu}}\:$
for e.q.a. because $\:e_{g1}e_{1h}\:=\:\frac{1}{2}(e_{gh}+e_{11}).\:$
Those families are minimal; certainly, all their extensions are also
parental sets for the same algebra.

{\bf Theorem 4.6.} {\em If $\:F\:=\:\{\:e_1,e_2\:\},\ e_1^2=e_1,\ 
e_2^2=e_2\:$ then the rank of the corresponding offspring subalgebra does 
not exceed 3.}

To be prepared for the proof below we consider a 
special nonnegative projection associated with $\:e_1\:$ by writing
$$
Bx\:=\:2e_1x\:-\:(2e_1x,e_1)e_1
\eqno (4.6)
$$
where $\:(,)\:$ is the standard inner product at the canonical basis
$\:\{\:e_i\:\}_1^n.\:$
This operator $\:B\:$ was introduced in [15] (see also [17], Sections 5.3 and
5.4). 
Let us recall its  properties, basically, without proofs. 
First of all, by (4.6)
$$
Be_1\:=\:0,\ \ {\rm Im}B\perp e_1.
\eqno (4.7)
$$

{\bf Lemma 4.7.} {\em In the algebra $\:\CA\:$ of type $\:(m,\gd)\:$ the
operator $\:B\:$ is a nonnegative projection of the form }
$$ 
Bx\:=\:\sum\limits_{i=1}^{m-1}(x,b_i^*)b_i
\eqno (4.8)
$$
{\em where}
$$
b_i\:\geq\:0,\ \ b_i^*\:\geq\:0,\ \ (b_i,b_k^*)\:=\:\gd_{ik}
\eqno (4.9)
$$
{\em and}
$$
{\rm supp}b_k\:\not\subset\:\bigcup\limits_{i\neq k}{\rm supp}b_i\ \ \ 
(1\leq k\leq m-1).
\eqno (4.10)
$$

{\bf Corollary 4.8.} rk$B\:=\:m-1.\:$

{\bf Proof.} It follows from (4.8) that Im$B\:\subset\:{\rm Lin}\{\:b_i\:\}_
1^{m-1}\:$ and actually these two subspaces coincide since 
$$
Bb_k\:=\:\sum\limits_{i=1}^{m-1}(b_k,b_i^*)b_i\:=\:b_k\ \ \ 
(1\leq k\leq m-1).
$$
Because of (4.9), $\:b_i^*\:$ are linearly indepedent. Thus, they form a 
basis in Im$B.\:$ 
$\ \ \ \ \ \ \ \ \ \ \ \ \ \ \ \ \ \Box\ $

>From now on we assume that $\:b_i\:$ are {\em normalized} in the sense 
$\:s(b_i)\:=\:1\ \ (1\leq i\leq m-1).\:$ Then we have

{\bf Lemma 4.9.} {\em The intersection $\:\gD_B\:=\:{\rm 
Im}B\:\bigcap\:\gD^{n-1}\:$ coincides with the convex hull of 
$\:\{\:b_i\:\}_1^{m-1}.\:$ }

Thus, $\:\gD_B\:$ is a symplex.

{\bf Proof.} Obviously, all convex combinations of $\:\{\:b_i\:\}_1^{m-1}\:$ 
belong to $\:\gD_B\:$ since all $\:b_i\:\in\:\gD_B.\:$
Conversely, if $\:x\:\in\:\gD_B\:$ then
$$
x\:=\:Bx\:=\:\sum\limits_{i=1}^{m-1}\ga_ib_i
$$
with $\:\ga_i\:=\:(x,b_i^*)\:\geq\:0\:$ and 
$\:\sum\limits\ga_i\:=\:s(x)\:=\:1.\:$
$\ \ \ \ \ \ \ \ \ \ \ \ \ \ \ \ \ \ \ \ \ \ \ \ \ \ \ \ \ \ \ \ \ \ \ \  
\ \ \ \ \ \ \ \ \ \ \ \ \ \ \ \ \ \ \ \ \ \Box\ $

Finally, we have

{\bf Lemma 4.10.} {\em The formulas}
$$
e_1b_i\:=\:\frac{e_1\:+\:b_i}{2}\ \ \ (1\leq i\leq m-1)
\eqno (4.11)
$$
{\em and}
$$
B(b_ib_k)\:=\:\frac{b_i\:+\:b_k}{2}\ \ \ (1\leq i,k\leq m-1)
\eqno (4.12)
$$
{\em hold. }

{\bf Proof of Theorem 4.6.} Without loss of generality we can assume 
that Lin$F^{'}\:=\:\CA\:$ and we have to prove that $\:m\:\leq\:3.\:$ 
Thus,
$$
e_1e_2\:=\:\sum\limits_{k=1}^n\pi_ke_k
\eqno (4.13)
$$     
where $\:\sum\limits\pi_k\:=\:1,\ \ \pi_k\geq 0\:$ and, moreover, 
$\:\pi_k >0\:$ for $\:k\geq 3\:$ since $\:\CA\:$ is the offspring algebra 
of $\:F\:=\:\{\:e_1,e_2\:\}.\:$ In terms of the projection $\:B\:$
$$
Be_2\:=\:2\sum\limits_{k=2}^n\pi_ke_k.
\eqno (4.14)
$$
On the other hand
$$
Be_2\:=\:2\sum\limits_{i=1}^{m-1}\gl_ib_i
\eqno (4.15)
$$
with $\:\gl_i=\frac{1}{2}(Be_2,b_i^*)\:\geq\:0,\ \ 1\leq i\leq m-1.\:$

Being a projection, $\:B=B^2,\:$ so (4.14) yields
$$
Be_2\:=\:2\sum\limits_{k=2}^n\pi_kBe_k
$$
therefore
$$
\gl_i\:=\:\sum\limits_{k=2}^n\pi_k(Be_k,b_i^*)\:=\:2\pi_2\gl_i\:+\:
\sum\limits_{k=3}^n\pi_k(Be_k,b_i^*)
\eqno (4.16)
$$
which implies $\:\gl_i >0\ \ (1\leq i\leq m-1).\:$ Indeed, if there is 
$\:\gl_i =0\:$
then $\:(Be_2,b_i^*)\:=\:0\:$ and 
$\:(Be_k,b_i^*)\:=\:0\ \ (3\leq k\leq n)\:$ 
from (4.16) and, finally, $\:(Be_1,b_i^*)\:=\:0\:$ from (4.7). As a result,
$\:(Bb_i,b_i^*)\:=\:0\:$ while $\:(Bb_i,b_i^*)\:=\:(b_i,b_i^*)\:=\:1.\:$

Now it follows from (4.15) that
$$
{\rm supp}(Be_2)\:=\:\bigcup\limits_{i=1}^{m-1}{\rm supp}b_i.
\eqno (4.17)
$$

Coming back to (4.13) and (4.14) we get
$$
e_1e_2\:=\:\pi_1e_1\:+\:\frac{1}{2}Be_2
\eqno (4.18)
$$
whence,
$$
(e_1e_2)^2\:=\:\pi_1^2e_1\:+\:\pi_1e_1(Be_2)\:+\:\frac{1}{4}(Be_2)^2.
$$
Multiply (4.15) by $\:e_1\:$ and using (4.11) we obtain
$$
e_1(Be_2)\:=\:e_1\sum\limits_{i=1}^{m-1}\gl_i\:+\:\sum\limits_{i=1}^{m-1}
\gl_ib_i\:=\:e_1\sum\limits_{i=1}^{m-1}\gl_i\:+\:\frac{1}{2}Be_2
$$
However, since $\:s(b_i)=1\ \ (1\leq i\leq m-1)\:$ we have
$$
\sum\limits_{i=1}^{m-1}\gl_i\:=\:\frac{1}{2}s(Be_2)\:=\:\sum\limits_{k=2}
^{m-1}\pi_k\:=\:1\:-\:\pi_1
$$
Finally,
$$
e_1(Be_2)\:=\:(1\:-\:\pi_1)e_1\:+\:\frac{1}{2}Be_2
$$
and then
$$
(e_1e_2)^2\:=\:\pi_1e_1\:+\:\frac{1}{2}\pi_1Be_2\:+\:\frac{1}{4}(Be_2)^2.
\eqno (4.19)
$$
On the other hand,
$$
(e_1e_2)^2\:=\:\frac{1}{4}e_1\:+\:\frac{1}{4}e_2\:+\:\frac{1}{2}e_1e_2
$$
by Proposition 2.2. With (4.18) this yields
$$
(e_1e_2)^2\:=\:(\frac{1}{4}\:+\:\frac{1}{2}\pi_1)e_1\:+\:\frac{1}{4}e_2\:+\:
\frac{1}{4}Be_2.
\eqno (4.20)
$$

Let us compare the $\:e_1$-coordinates in (4.19) and (4.20). Taking into 
account 
that $\:Be_2\:\perp\:e_1\:$ (see (4.7)) and $\:(Be_2)^2\:\geq\:0\:$  
we get $\:\frac{1}{4}\:+\:\frac{1}{2}\pi_1\:\geq\:\pi_1,\:$
i.e. $\:\pi_1\leq\frac{1}{2}.\:$ Now we compare the $\:e_2$-coordinates
and get
$$
\frac{1}{4}\:+\:\frac{1}{2}\pi_2\:=\:\pi_1\pi_2\:+\:
\frac{1}{4}((Be_2)^2,e_2)
$$
whence,
$$
\frac{1}{4}\:+\:(\frac{1}{2}-\pi_1)\pi_2\:=\:\frac{1}{4}
((Be_2)^2,e_2)\:=\:\sum\limits_{i,k=1}^{m-1}\gb_{ik}\gl_i\gl_k
\eqno (4.21)
$$
where $\:\gb_{ik}\:=\:(b_ib_k,e_2)\:$ by (4.15). Since all $\:b_i\geq0\:$ 
we have $\:\gb_{ik}\geq0.\:$ Moreover, there exists $\:\gb_{i_1k_1}>0,\:$
otherwise all $\:\gb_{ik}=0\:$ which contradicts (4.21) because of
$\:\frac{1}{2}-\pi_1\geq0\:$ and $\:\pi_2\:\geq\:0.\:$ 
Applying $\:B\:$ to the inequality $\:b_{i_1}
b_{k_1}\geq\gb_{i_1k_1}e_2\:$ and using (4.12) we obtain
$$
\frac{b_{i_1}\:+\:b_{k_1}}{2}\:\geq\:\gb_{i_1k_1}Be_2.
$$
Hence,
$$
{\rm supp}(Be_2)\:\subset\:{\rm supp}b_{i_1}\:\bigcup\:{\rm supp}b_{k_1}.
$$
However, we have (4.17). Therefore
$$
\bigcup\limits_{i=1}^{m-1}{\rm supp}b_i\:=\:{\rm supp}b_{i_1}\:\bigcup\:
{\rm supp}b_{k_1}.
$$
In view of (4.10) we get $\:m=2\:$ if $\:i_1=k_1\:$ and $\:m=3\:$ if
$\:i_1\neq k_1.\:$
$\ \ \ \ \ \ \ \ \ \ \ \ \ \ \ \ \ \ \ \ \ \ \ \ \ \ \ \ \ \ \ \Box\ $

As a consequence we obtain a very useful 

{\bf Theorem 4.11.} {\em Let $\:F\:=\:\{\:e_1,e_2\:\},\ e_1^2=e_1,\ 
e_2^2=e_2.\:$ If the correspondidg offspring subalgebra is normal then it is
either 2-dimensional u.a. or 3-dimensional e.M.a., or g.a. (which is 
4-dimensional).}

{\bf Proof.} By Theorem 4.6 the subalgebra is of rank $\leq\:3.\:$ By 
Theorem 2.11 this is regular. By Theorem 3.1 this is either e.M.a. or e.q.a.

If a normal e.M.a. is the offspring subalgebra of two parental 
idempotents then its dimension is 2 or 3. This is an u.a. in the case of 
dimension 2. If an e.q.a. has rank $\leq\:3\:$ then it is q.a.
$\ \ \ \ \ \ \ \ \ \ \ \ \ \ \ \ \ \ \ \ \ \ \ \ \ \ \ \ \ \ \ \ \ \ \
\ \ \ \ \ \ \ \ \ \ \ \ \ \ \ \ \ \ \ \ \ \ \ \ \ \ \ \ \ \ \ \ \ \ \ \ \ 
\ \ \ \ \ \ \ \ \ \ \ \ \ \ \ \ \ \ \ \ \ \ \ \ \ \ \Box\ $

We say that a stochastic Bernstein algebra is {\em grounded} if it is the 
offspring subalgebra of the set of its basis idempotents. As we know 
every normal stochastic regular algebra is grounded.

{\bf Lemma 4.12.} {\em Let $\:\CA\:$ be not constant. If every proper 
coordinate subalgebra of $\:\CA\:$ is grounded then the greatest 
offspring subalgebra $\:\CA_0\:$ is also grounded.}

Obviously, $\:\CA_0\:$ is the linear span of the set
$$
\bigcup\limits_{1\leq i\leq k\leq n}{\rm supp}(e_ie_k).
\eqno (4.22)
$$
The corresponding invariant face $\:\gD\:$ of the simplex $\:\gD^{n-1}\:$
is the convex hull of the same set (4.22). On the other hand,$\:\gD\:$
is the smallest face containing Im$V,\:$ the image of the evolutionary 
operaror $\:Vx\:=\:x^2\ \ (x\:\in\:\gD^{n-1}).\:$ Indeed, 
$$
x^2\:=\:\sum\limits_{i,k=1}^nx_ix_ke_ie_k
$$
for
$$
x\:=\:\sum\limits_{i=1}^nx_ie_i.
$$
Therefore
$\:x\:\in\:{\rm Im}V\:\Rightarrow\:x=x^2\:\in\:\gD\:$ and 
$\:x\:\in\:{\rm Int}\gD^{n-1}\:\Rightarrow\:x^2\:\in\:{\rm Im}V\:\bigcap\:
{\rm Int}\gD,\:$
so that Im$V\:\subset\:\gD\:$ and Im$V\:\bigcap\:{\rm Int}\gD\:\neq\:
\emptyset.\:$

A face $\:\gG\:$ of the simplex $\:\gD^{n-1}\:$ is called {\em essential} 
if $\:C_{\gG}\:\equiv\:{\rm Im}V\:\bigcap\:{\rm Int}\gG\:\neq\:\emptyset.\:$
Obviously,
$$
{\rm Im}V\:=\:\bigcup\limits_{\gG}C_{\gG}\:=\:\bigcup\:\{\:C_{\gG}\: 
:\:\gG\  {\rm is\  essential}\:\}.
\eqno (4.23)
$$
It turns out that this partition is an {\em elementary cell complex on}
Im$V\:$ in the following sense (see [14]; [17], Section 5.7).

For any Hausdorff topological space $\:X\:$ a subset $\:C\:\subset\:X\:$ is 
called a $\:\nu$-{\em dimensional elementary cell} if there exists a 
bounded open set $\:U\:\subset\:\R^{\nu}\:$ whose closure $\:\overline{U}\:$
is contractible (within itself to a point) and homeomorphic to 
$\:\overline{C}\:$ by a boundary preserved homeomorphism.

A finite partition of $\:X\:$ is called an {\em elementary cell complex on}
$\:X\:$ if 1) all the parts are elementary cells; 2) the boundary of each
of one is a union of some lower dimensional cells; 3) the intersection of 
the closures of any two cells is contractible.

The maximal cell dimension $\:d\:$ is called the {\em dimension of the 
complex.} (It is equal to the usual topological dim$X).\:$

In our case $\:X\:=\:{\rm Im}V,\:$ the cells are $\:C_{\gG}\:$ for essential
faces $\:\gG\:$ and dim$C_{\gG}\:=\:m_{\gG}\:-\:1\:$ where $\:m_{\gG}\:=\:
{\rm rk}(\CA_{V\mid\gG}).\:$ If $\:\gG_1\:\subset\:\gG_2\:$ and $\:\gG_1\:
\neq\:\gG_2\:$ then dim$C_{\gG_1}\: <\:{\rm dim}C_{\gG_2}.\:$ Thus, the
{\em dimension of this elementary cell complex is} $\:m-1\:$ where 
$\:m\:=\:{\rm rk}\CA\:$ as usual. The only $\:(m-1)$-dimensional cell is
$\:C_{\gD}\:=\:{\rm Im}V\:\bigcap\:{\rm Int}\gD.\:$

{\em Every essential face $\:\gG\:$ is invariant} and
$$
{\rm Im}(V\mid\gG)\:=\:{\rm Im}V\:\bigcap\:\gG\:=\:\overline{C_{\gG}}.
\eqno (4.24)
$$
The boundary $\:\partial 
C_{\gG}\:=\:\overline{C_{\gG}}\:\setminus\:C_{\gG}\:$ is actually
$$
\partial C_{\gG}\:=\:{\rm Im}V\:\bigcap\:\partial \gG.
\eqno (4.25)
$$
Obviously, $\:\partial C_{\gG}\:=\:\emptyset\:$ iff dim$C_{\gG}\:=\:0\:$ 
which means that $\:V\mid\gG\:$ is constant.

For a topological reason ([17]), Lemma 5.7.2) in any $\:d$-dimensional 
elementary cell complex the number of 0-dimensional cells is at least
$\:d+1.\:$ Hence {\em there exists at least $\:m\:$ constant
coordinate subalgebras in any stochastic Bernstein algebra} ([14]; [17], 
Theorem 5.7.1)

{\bf Remark 4.13.} At least one constant subalgebra can be obtained in a 
much more simple way. This follows from Theorem 5.2.1 [17] saying that {\em 
the subalgebra corresponding to a minimal invariant face is constant.}
(The latter is a generalization of a Bernstein theorem provided in [4] 
with a very complicated proof. A short proof was found in [9].)

After these preliminaires we can directly pass to

{\bf Proof of Lemma 4.12.} As aforesaid, there exists a constant 
coordinate subalgebra in $\:\CA.\:$ This is a proper subalgebra because 
$\:\CA\:$ is not constant. Being grounded this subalgebra is 1-dimensional,
generated by a basis idempotent. We conclude that the set of all basis 
idempotents is not empty. Let $\:\CA_1\:$ be its offspring subalgebra, so
$\:\CA_1\:$ is the greatest grounded subalgebra. Obviously, $\:\CA_1\:
\subset\:\CA_0.\:$ We have to prove that $\:\CA_1\:=\:\CA_0.\:$

Suppose that $\:\CA_1\:\neq\:\CA_0.\:$ Then $\:\gG_1\:\neq\:\gD\:$ where
$\:\gG_1\:$ is the invariant face corresponding to the algebra $\:\CA_1.\:$
Since $\:\gG_1\:\subset\:\gD\:$ and $\:\gG_1\:\neq\:\gD,\:$ we have
$\:\gG_1\:\subset\:\partial \gD.\:$ By (4.24) and (4.25)
$$
\overline{C_{\gG_1}}\:=\:{\rm Im}V\:\bigcap\:\gG_1\:\subset\:{\rm Im}V\:
\bigcap\:\partial\gD\:=\:\partial C_{\gD}.
\eqno (4.26)
$$
However, $\:\partial\gD\:=\:\bigcap{\rm Int}\gG\:$ where $\:\gG\:$ runs
over all faces $\:\gG\:\subset\:\gD,\ \ \gD\:\neq\:\gD^{n-1}.\:$
Hence,
$$
\partial C_{\gD}\:=\:\bigcup\:\{\:{\rm Im}V\:\bigcap\:{\rm Int}\gG\: :\:
\gG\:\subset\:\gD, \ \ \gG\:\neq\:\gD\:\}\:=\:\bigcup\:\{\:C_{\gG}\: :\:
\gG\:\subset\:\gD,\ \ \gG\:\neq\:\gD\:\}.
$$
Since $\:\partial C_{\gD}\:$ is closed, we get
$$
\partial C_{\gD}\:=\:\bigcup\:\{\{\overline{C_{\gG}}\: :\:\gG\:\subset\:\gD,
\ \gG\:\neq\:\gD\:\}.
\eqno (4.27)
$$
All $\:\gG\:$ in (4.27) may be supposed to be essential (otherwise 
$\:C_{\gG}\:=\:\emptyset\:$). Therefore they are invariant, i.e. they 
correspond 
to some coordinate subalgebras. Being proper these subalgebras are grounded
hence, they are contained in $\:\CA_1.\:$ Hence, $\:\gG\:\subset\:\gG_1\:$
for all essential $\:\gG\:$ in (4.27). Respectively, $\:\overline{C_{\gG}}\:
\subset\:\overline{C_{\gG_1}}\:$ and we conclude that $\:\partial C_{\gD}\:
\subset\:\overline{C_{\gG_1}}.\:$ Jointly with (4.26) this results in the
equality $\:\partial C_{\gD}\:=\:\overline{C_{\gG_1}}.\:$ But this 
contradicts a well known topological fact: the boundary of any cell (except
for 0-dimensional one) is not contractible.
$\ \ \ \ \ \ \ \ \ \ \ \ \ \ \ \ \ \ \ \ \ \ \ \ \ \Box\ $

{\bf Corollary 4.14.} {\em Let $\:\CA\:$ be not constant and nondegenerate.
If every proper coordinate subalgebra of $\:\CA\:$ is grounded then 
$\:\CA\:$ is also grounded.}

{\bf Proof.} The nondegeneracy means that $\:\gD\:=\:\gD^{n-1},\:$ i.e.
$\:\CA_0\:=\:\CA.\:$
$\ \ \ \ \ \ \ \ \ \ \ \ \ \ \ \ \ \ \ \ \ \ \ \ \ \Box\ $

\section{Proof of the Main Theorem}
$\ \ \ \ \ \ \ \ $
{\em Given an ultranormal stochastic Bernstein algebra $\:\CA.\:$ We have
to prove that $\:\CA\:$ is regular.} 
As usual, $\:(m,\gd)\:$ denotes the type of $\:\CA,\:$ 
dim$\CA\:=\:n\:=\:m+\gd.\:$

Above all, let us come back to the projection $\:B\:$ which is associated 
with a basis idempotent, say $\:e_1,\:$ via (4.6).
Such an idempotent does exist because all constant subalgebras of 
$\:\CA\:$ are 1-dimensional by ultranormality. As we know, a constant 
subalgebra does exist (Remark 4.13), moreover, there exist at least $\:m\:$
constant subalgebras, so there are at least $\:m\:$ basis idempotents in
$\:\CA.\:$ In fact, some $\:m\:$ basis idempotents can be obtained by one 
of them using the projection $\:B.\:$ This way also yiels an additional
useful information.

{\bf Lemma 5.1.} {\em In notation of Lemma 4.7, for every vector $\:b_i\ 
\ (1\leq i\leq m-1)\:$ there exists a unique basis idempotent 
$\:e_{j_i},\ \ j_i>1,\:$ such that $\:Be_{j_i}\:=\:\gl_ib_i,\ \ 
\gl_i\:>\:0.\:$}

{\bf Proof.} Since $\:b_i^2\:$ is a nonzero idempotent and 
$\:b_i^2\:\geq\:0,\:$ the subspace $\:L\:=\:{\rm Lin(supp}b_i^2)\:$ is a
subalgebra. Being a coordinate subalgebra of the ultranormal algebra
$\:\CA,\ L\:$ is normal. If it is constant then dim$L\:=\:1\:$ i.e. 
$\:b_i^2\:=\:
e_j\:$ where $\:e_j\:$ is a basis vector (recall that $\:s(b_i)=1),\:$ in
fact, $\:e_j\:$ is a basis idempotent. By (4.12) $\:Be_j\:=\:b_i.\:$

Let $\:L\:$ be nonconstant. Then rk$L\geq2\:$ hence, $\:L\:$ has at least 
two of basis idempotents. One of them is not $\:e_1,\:$ 
say it is $\:e_2,\:$ so that $\:e_2\:\in{\rm supp}(b_i^2).\:$ This means that
$\:\gb_2\:>\>0\:$ in the expansion
$$
b_i^2\:=\:\sum\limits_{j=1}^n\gb_je_j.
$$
Applying $\:B\:$ we get by (4.12)
$$
b_i\:=\:\sum\limits_{j=2}^n\gb_jBe_j.
$$
However, $\:b_i\:$ is an extreme point in the symplex $\:\gD_B\:$
(see Lemma 4.9). Hence $\:Be_2\:=\:\gl b_i\:$ with $\:\gl\geq 0.\:$ Actually
$\:\gl >0\:$ because $\:Be_2\:=\:0\:$ means that $\:e_1e_2\:=\:e_1\:$ and 
then
Lin$\:\{\:e_1,e_2\:\}\:$ is a subalgebra but not unit and nonconstant which
is impossible.
  
It remains to prove that $\:e_2\:$ is the only idempotent such that
$\:Be_2\:=\:\gl b_i,\ \gl >0.\:$

The last equality means that $\:2e_1e_2\:=\:\ga e_1\:+\:\gl b_i,\ \ \ga =
1-\gl.\:$ But, according to Theorem 4.6 we have only three cases:
1) $\:2e_1e_2\:=\:e_1\:+\:e_2\:$ (u.a.); 2) $\:2e_1e_2\:=\:\ga e_1\:+\:
\gb e_2\:+\:\gamma e_3,\ \ \gamma >0\:$ (e.M.a.); 3) $\:2e_1e_2\:=\:e_3\:+\:
e_4\:$ (q.a.). Respectively, $\:Be_2=e_2\:$ or $\:Be_2=\gb e_2+\gamma 
e_3\ (\gamma >0),\:$ or $\:Be_2= e_3+e_4,\:$ i.e. $\:\gl b_i=e_2\:$ or
$\:\gl b_i=\gb e_2+\gamma e_3\ (\gamma >0),\:$ or $\:\gl b_i=e_3+e_4.\:$

Since $\:s(b_i)=1\:$ we get such three cases: 1) $\:b_i=e_2;\:$ 2) $\:b_i=
\gep e_2+\gw e_3\ (\gw >0);\:$ 3) $\:b_i=\frac{1}{2}(e_3+e_4).\:$ In case
1) $\:b_i\:$ coincides with the idempotent $\:e_2.\:$ In case 2)
$\:e_3\:$ is the only nonidempotent in supp$b_i\:$ and
$\:e_2\:$ is the only idempotent in supp$e_3^2\:$ different from $\:e_1\:$
(see (3.11)). Finally, in case 3) supp$b_i\:$ consists of two 
idempotents, $\:e_3\:$ and $\:e_4\:$ and 
$\:e_2\:=\:4b_i^2\:-\:2b_i\:-\:e_1.\:$
$\ \ \ \ \ \ \ \ \ \ \ \ \ \ \ \ \ \ \ \ \ \ \ \ \ \ \ \ \ \ \
\ \ \ \ \ \ \ \ \ \ \ \ \ \ \ \ \ \ \ \ \ \ \ \ \ \ \ \ \ \ \ 
\ \ \ \ \ \ \ \ \ \ \ \ \ \ \ \ \ \ \ \ \ \ \ \ \ \Box\ $

It is convenient to denote the algebras in case 1), 2) and 3) by
$\:\{\:e_1,e_2\mid\emptyset\:\},\ \{\:e_1,e_2\mid e_3\:\}\:$ and $\:\{\:
e_1,e_2\mid e_3,e_4\:\}\:$ respectively.

Let the number of idempotents in the canonical basis is $\:\rho,\:$  
so we can assume that they are $\:e_1,\:...,\:e_{\rho}.\:$ We already
know that $\:\rho\geq m\:$ (following Lemma 5.1 or the previous 
topological argumentation). Now we even get

{\bf Corollary 5.2.} rk$\:\{\:Be_i\:\}^{\rho}_2\:=\:m-1.\:$

{\bf Proof.} As we know $\:\{\:b_i\:\}_1^{m-1}\:$ is a basis in Im$B.\:$ On
the other hand, $\:\{\:b_i\:\}_1^{m-1}\:\subset\:{\rm 
Lin}\:\{\:Be_i\:\}_2^{\rho}\:$ by Lemma 5.1.
$\ \ \ \ \ \ \ \ \ \ \ \ \ \ \ \ \ \ \ \ \ \ \ \ \ \ \ \ \ \ \ \ \ \ \ \ 
\ \ \ \ \ \ \ \ \ \ \ \ \ \ \ \ \ \ \ \ \
\ \ \ \ \ \ \ \ \ \ \ \ \ \ \ \ \ \ \ \ \ \ \ \ \Box\ $

We continue the proof of the Main Theorem in the frameworks of the 
following alternative: {\em all basis vectors $\:e_i\ (1\leq i\leq n)\:$ 
are idempotents, i.e. $\:\rho =n,\:$ or not, i.e. $\:\rho <n.\:$} Let us say
that the first possibility is the {\em pure case} and the second one is 
the {\em mixed case.}

{\bf A). The pure case $\:(\rho =n).\:$}  
Take a pair $\:\{\:e_i,e_k\:\}\:$ of basis 
vectors, $\:i\neq k.\:$ Its offspring subalgebra is normal as every 
coordinate subalgebra of $\:\CA.\:$ According to Theorem 4.11 this is either
$\:\{\:e_i,e_k\mid\emptyset\:\}\:$ or $\:\{\:e_i,e_k\mid e_g,e_h\:\}.\:$ 
Using the multiplication $\:R\:$ introduced by (2.17) we 
have $\:R(e_i,e_k)\:=\:0\:$ or $\:R(e_i,e_k)\:\neq\:0\:$ respectively and 
in the second case
$\:R(e_i,e_g)\:=\:0,\ \ R(e_i,e_h)\:=\:0\:$ and similarly for $\:e_k.\:$

{\bf Lemma 5.3} {\em If $\:R(e_i,e_k)\:\neq\:0\:$ then there is no $\:e_j\:$
such that $\:R(e_i,e_j)\:=\:0\:$ and $\:R(e_k,e_j)\:=\:0\:$ except for 
$\:e_j=e_g\:$ and $\:e_j=e_h.\:$ } 

{\bf Proof.} The subspace $\:L\:=\:{\rm Lin}\:\{\:e_j,e_i,e_k,e_g,e_h\:\}\:$
is the offspring subalgebra of the family $\:\{\:e_j,e_i,e_k\:\}.\:$
The algebra $\:L\:$ is nuclear as the linear span of a set of idempotents. By
Corollary 2.12 $\:L\:$ is regular. Moreover, $\:L\:$ is normal. By 
Theorem 3.1 $\:L\:$ must be an e.q.a. Indeed, $\:L\:$ is not an u.a. since
$\:L\:$ contains a q.a. and $\:L\:$ is not an e.M.a. since all vectors from
the canonical basis of $\:L\:$ are idempotents. Under $\:4\leq{\rm dim}L\leq
5,\:$ actually dim$L=4\:$ because no prime number can be dimension of an 
e.q.a. Since $\:e_j\neq e_i\:$ and $\:e_j\neq e_k,\:$ we conclude that
$\:e_j =e_g\:$ or $\:e_j =e_h.\:$
$\ \ \ \ \ \ \ \ \ \ \ \ \ \ \ \ \ \ \ \ \ \ \ \ \ \ \ \ \ \ \ \ \ \ \ \ \
\ \ \ \ \ \ \ \ \ \ \ \ \ \ \ \ \ \ \ \ \ \ \ \ \ \ \ \ \ \ \ \ \ \
\ \ \ \ \ \ \ \ \ \ \ \ \ \ \ \ \ \ \ \ \ \ \ \ \Box\ $

Let us write $\:e_iR_0e_k\:$ in the case $\:R(e_i,e_k)\:=\:0,\:$
so that $\:R_0\:$ is a binary relation on the set 
$\:\{\:e_j\:\}_1^n.\:$
Obviosly, it is reflexive and symmetric. For any $\:e_j\:$
we define its {\rm pool}
$$
P(e_j)\:=\:\{\:e_k\: :\:e_jR_0e_k\:\}\:=\:\{\:e_k\: 
:\:R(e_j,e_k)=0 \:\}
\eqno (5.1)
$$
We also consider the {\em punctured pool}
$$
P^*(e_j)\:=\:P(e_j)\setminus\{\:e_j\:\}.
$$

{\bf Lemma 5.4.} {\em The equality }
$$
{\rm card}P^*(e_j)\:=\:m-1\ \ \ (1\leq j\leq n)
\eqno (5.3)
$$
{\em holds.}

{\bf Proof.} The projection $\:B_j\:$ associated with $\:e_j\ \ (B_1\:=\:B\:$
in this notation) acts as follows:
$$
B_je_k\:=\:e_k\ \ \ (e_k\:\in\:P^*(e_j))
\eqno (5.4)
$$
and
$$
B_je_k\:=\:e_g\:+\:e_h
\eqno (5.5)
$$
if $\:e_k\:\not\in\:P^*(e_j)\:$ and $\:\{\:e_k,e_j\mid e_g,e_h\:\}\:$ 
is the corresponding q.a. Since $\:e_g\:$ and $\:e_h\:$ belong to
$\:P^*(e_j),\:$ it follows from (5.4) and (5.5) that rk$B_j\:=\:
{\rm card}P^*(e_j).\:$ On the other hand, rk$B_j=m-1.\:$
$\ \ \ \ \ \ \ \ \ \ \ \ \ \ \ \ \ \ \ \ \ \ \ \ \ \ \ \ \ \ \ \ 
\ \ \ \ \ \ \ \ \ \ \ \ \ \ \ \ \ \ \ \ \ \ \ \ \ \ \ \ \ \ \ \ \
\ \ \ \ \ \ \ \ \ \ \ \ \ \ \ \ \ \ \ \ \ \ \ \ \ \ \ \ \ \ \ \ \Box $

{\bf Corollary 5.5.}  card$P^*(e_j)\:$ {\em is independent of $\:j.\:$ }

Coming back to the binary relation $\:R_0\:$ we prove

{\bf Lemma 5.6.} {\em The restriction $\:R_0\mid P^*(e_j)\:$ is
an equivalence relation.}

{\bf Proof.} We only need to check that $\:R_0\:$ is transitive on $\:
P^*(e_j).\:$ Let $\:\{\:e_i,e_k,e_l\:\}\:$ be a triple from $\:P^*(e_j)\:$
such that $\:R(e_i,e_l)\:=\:0\:$ and $\:R(e_k,e_l)\:=\:0\:$ but $\:R(e_i,e_k)
\:\neq\:0.\:$ Since $\:R(e_i,e_j)\:=\:0\:$ and $\:R(e_k,e_j)\:=\:0\:$ as 
well, Lemma 5.3 yields the q.a. $\:\{\:e_i,e_k\mid e_j,e_l\:\}\:$ which 
contradicts $\:R(e_j,e_l)\:=\:0.\:$
$\ \ \ \ \ \ \ \ \ \ \ \ \ \ \ \ \ \ \ \ \ \ \ \ \ \ \ \ \ \ \ \ \ \ \ 
\ \ \ \ \ \ \ \ \ \ \ \ \ \ \ \ \ \ \ \ \ \ \  
\ \ \ \ \ \ \ \ \ \ \ \ \ \ \ \ \ \ \ \ \ \ \ \ \Box\ $

Obviously, all classes of this equivalence relation are  u.a.

From now on we {\em assume that $\:\CA\:$ is not unit.} (Otherwise, $\:\CA\:$
is regular a fortiori.)

{\bf Lemma 5.7.} {\em There are exactly two classes of the relation $\:R_0
\mid P^*(e_j),\ \ 1\leq j\leq n.\:$}

(By the way, we see that $\:R_0\:$ is not transitive on the 
whole pool $\:P (e_j)).\:$

{\bf Proof.} Let $\:e_l\:\not\in\:P (e_j),\:$ so that $\:R(e_j,e_l)\:\neq
\:0.\:$ Then we have $\:\{\:e_j,e_l\mid e_i,e_k\:\}\:$ where $\:e_i\:$ 
and $\:e_k\:$ belong to $\:P^*(e_j)\:$ and $\:R(e_i,e_k)\:\neq\:0,\:$
so that $\:e_i\:$ and $\:e_k\:$ are not equivalent. In such a way either 
there are at least two classes in $\:P^*(e_j)\:$ or $\:P (e_j)=\CA\:$
and $\:P^*(e_j)\:$ is an entire class. But in the last case $\:P^*
(e_j)\:$ an u.a. and then $\:\CA\:$ is so as.

Suppose that there are more than two classes in $\:P^*(e_j).\:$ Then 
there is a triple $\:\{\:e_i,e_k,e_l\:\}\:\subset\:P^*(e_j)\:$ such that
$\:R(e_i,e_k)\:\neq\:0,\ \ R(e_i,e_l)\:\neq\:0\  {\rm and}\ 
R(e_k,e_l)\:\neq\:0.\:$ By Lemma 5.3 there are three q.a., namely,
$$
\{\:e_i,e_k\mid e_j,e_p\:\},\ \ \{\:e_i,e_l\mid e_j,e_q\:\},\ \ 
\{\:e_k,e_l\mid e_j,e_r\:\}.
\eqno (5.6)
$$
In (5.6) the seven involved vectors are pairwise distinct. For example, 
$\:e_p\neq e_q\:$
since $\:2e_pe_j\:=\:e_i\:+\:e_k\ {\rm but}\ 2e_qe_j\:=\:e_i\:+\:e_l.\:$
Also $\:e_p\neq e_l\ {\rm since}\ R(e_k,e_l)\:\neq\:0\ {\rm but}\ R(e_k,e_p)
\:=\:0.\:$
In situation (5.6) the punctured pool $\:P^*(e_j)\ {\rm is}\ 
\{\:e_i,e_k,e_l\:\}.\:$ According to Lemma 5.4 $\:m=4,\:$ so card 
$P^*(e_i)\:=\:3\:$ by Corollary 5.5. However, 
$\:P^*(e_i)\:\supset\:\{\:e_j,e_p,e_q\:\}\:$ hence,
$$
P^*(e_i)\:=\:\{\:e_j,e_p,e_q\:\}.
\eqno (5.7)
$$
Since $\:e_r\:\not\in\:P^*(e_i),\ {\rm i.e.}\ R(e_i,e_r)\:\neq\:0,\:$ 
we get one more q.a., say $\:\{\:e_i,e_r\mid e_g,e_h\:\},\:$ where
$$
\{\:e_g,e_h\:\}\:\subset\:P^*(e_i)\:\bigcap\:P^*(e_r).
\eqno (5.8)
$$
This is a contradiction because $\:P^*(e_r)\:\supset\:\{\:e_k,e_l\:\}\:$
(see (5.6)) and card$\:P^*(e_r)\:=\:3\:$ so the intersection (5.8) 
can not contain more than one element.
$\ \ \ \ \ \ \ \ \ \ \ \ \ \ \  
\ \ \ \ \ \ \ \ \ \ \ \ \ \ \ \ \ \ \ \ \ \ \ \ \ \ \ \ \Box\ $

Let us denote the classes of $\:R_0\mid P^*(e_j)\ {\rm by}\ C_j\ 
{\rm and}\ \overline{C_j}.\:$

{\bf Lemma 5.8.} {\em There exists a bijective mapping from the complement
of the pool $\:P (e_j)\:$ onto the Cartesian product $\:C_j\times
\overline{C_j},\ \ 1\leq j\leq n.\:$}

{\bf Proof.} 
For any $\:e_k\:\not\in\:P(e_j)\:$ we have the q.a. $\:\{\:e_j,e_k\mid 
e_{g_k},e_{h_k}\:\}\:$ where $\:e_{g_k}\ {\rm and}\ e_{h_k}\:$ are both from
the punctured pool $\:P^*(e_j)\ {\rm and}\ R(e_{g_k},e_{h_k})\:\neq\:0\:$
which means that $\:e_{g_k}\ {\rm and}\ e_{h_k}\:$ are from different 
classes, say $\:e_{g_k}\:\in\:C_j\ {\rm and}\ 
e_{h_k}\:\in\:\overline{C_j.}\:$ The mapping defined in such a way is
injective since $\:e_k\:=\:2e_{g_k}e_{h_k}\:-\:e_j.\:$ It is also 
surjective. Indeed, if $\:e_g\:\in\:C_j\ {\rm and}\ e_h\:\in\:
\overline{C_j}\:$ then $\:R(e_g,e_h)\neq 0,\:$ so we have the q.a. 
$\:\{\:e_g,e_h\mid e_j,e_k\:\}\:$ where $\:e_j\:$ appears by Lemma 5.3
and then $\:e_k\:\not\in\:P (e_j).\:$ This means that $\:e_g\:=\:e_{g_k}
\ {\rm and}\ e_h\:=\:e_{h_k}.\:$
$\ \ \ \ \ \ \ \ \ \ \ \ \ \ \ \ \ \ \ \ \ \ \ \ \ \ \ \Box\ $

{\bf Corollary 5.9.} {\em Let $\:m_j\:=\:{\rm card}C_j+1\ {\em and}\ 
\overline{m_j}\:=\:{\rm card}\overline{C_j}+1.\:$ Then}
$$
m_j\:+\:\overline{m_j}\:=\:m+1,\ \ \ m_j\overline{m_j}\:=\:n\ \ \ 
(1\leq j\leq n).
\eqno (5.9)
$$

{\bf Proof.} Since $\:C_j\:\bigcup\:\overline{C_j}\:=\:P^*(e_j)\:$ and
$\:C_j\:\bigcap\:\overline{C_j}\:=\:\emptyset\:$ we obtain the first of 
equalities
(5.9) from (5.3). On the other hand, it follows from Lemma 5.8 that
$$
\ \ \ \ \ \ \ \ \ \ \
n\:=\:{\rm card}(C_j\:\times\:\overline{C_j})\:+\:{\rm card}P (e_j)\:=\:
(m_j-1)(\overline{m_j}-1)\:+\:m\:=\:m_j\overline{m_j}.
\ \ \ \ \ \ \ \ \ \ \ \ \ \ \ \Box
$$

{\bf Corollary 5.10.} $\:\gd\:=\:(m_j-1)(\overline{m_j}-1).\:$

Since $\:\CA\:$ is not unit, we have $\:\gd\:>\:0.\:$ Therefore 
$\:m_j\:\geq\:2\ {\rm and}\ \overline{m_j}\:\geq\:2.\:$

{\bf Corollary 5.11.} {\em The numbers $\:m_j\ {\em and}\ \overline{m_j}\:$
do not depend on $\:j.\:$}

Therefore one can set $\:m_j\:=\:\nu,\ \overline{m_j}\:=\:\overline{\nu}\:$
for all of $\:j,\ 1\leq j\leq n.\:$

Now we enumerate the basis $\:\{\:e_j\:\}\:$ in a new way:
$$
e_1\:\equiv\:e_{11},\ \ C_1\:=\:\{\:e_{i1}\:\}_{i=2}^{\nu},\ \ 
\overline{C_1}\:=\:\{\:e_{1k}\:\}_{k=2}^{\overline{\nu}}.
\eqno (5.10)
$$
By Lemma 5.8 the ordered pairs $\:(e_{i1},e_{1k})\:$ are in $\:1-1\:$
correspondence with the complement of the pool $\:P (e_1).\:$ Hence,
this complement can be listed as $\:\{\:e_{ik}\: :\:2\leq i\leq\nu,\ 
2\leq k\leq\overline{\nu}\:\}\:$ The correspondence is established by the
q.a. $\:\{\:e_{11},e_{ik}\mid e_{i1},e_{1k}\:\}.\:$ The whole basis 
$\:\{\:e_j\:\}_1^n\:$ can be written in the matrix form,
$\:E\:=\:(e_{ik})\ \ \ (1\leq i\leq\nu,\ \ 1\leq k\leq\overline{\nu}).\:$

{\bf Lemma 5.12.} {\em For any element $\:e_{ik}\:$ its pool $\:P 
(e_{ik})\:$ is the union of the i-th row and the k-th column of the matrix 
$\:E.\:$ Being punctured at $\:e_{ik}\:$ these lines are the equivalence 
classes of the punctured pool $\:P^*(e_{ik}).\:$ }

{\bf Proof.} Let us denote the punctured $\:k$-th column and $\:i$-th row by
$\:C_{ik}\:$ and $\:\overline{C_{ik}}\:$ respectively. In particular,
$\:C_{11}\:=\:C_1\ {\rm and}\ \overline{C_{11}}\:=\:\overline{C_1}\:$ by 
(5.10), so the lemma is true for $\:P (e_{11}).\:$ Now we consider
$\:P (e_{i1}),\ i>1.\:$

The 1-st column is $\:C_1\:\bigcup\:\{\:e_{11}\:\}\:$ therefore $\:R(e_{i1},
e_{g1})\:=\:0\:$ for all of $\:g,\ 1\leq g\leq\nu.\:$ Thus, $\:C_{i1}\:$ 
is contained in a class of $\:P^*(e_{i1}).\:$ Note that $\:e_{11}
\:\in\:C_{i1}.\:$ As the q.a. $\:\{\:e_{11},e_{ik}\mid e_{i1},e_{1k}\:\}\ 
\ (i,k\neq 1)\:$ shows, $\:R(e_{i1},e_{ik})\:=\:0\:$ for all of $\:k,\ 
1\leq k\leq\overline{\nu}.\:$ Thus, $\:\overline{C_{i1}}\:$ is contained
in $\:P^*(e_{i1}).\:$ However, $\:R(e_{11},e_{ik})\:\neq\:0\:$ for
$\:k\neq 1\:$ hence, $\:\overline{C_{i1}}\:$ lies in another class. In fact,
$\:C_{i1}\:$ and $\:\overline{C_{i1}}\:$ must coincide with the 
corresponding classes because of the same (up to transposition, a priori)
cardinalities. The lemma is proved for $\:P (e_{i1}).\:$
Quite similarly, this is true for $\:P (e_{1k}).\:$ But then $\:R(e_{ik},
e_{jk})=0\:$ for all $\:j,\ 1\leq j\leq\nu\:$ and $\:R(e_{ik},e_{il})=0\:$
for all $\:l,\ 1\leq l\leq\overline{\nu}\:$ which yields the first part 
of the lemma for $\:P(e_{ik}).\:$ Moreover, the second part is also
true. Indeed, $\:C_{ik}\:=\:(C_{ik}\:\bigcup\:\{\:e_{1k}\:\})\:\setminus\:
\{\:e_{ik}\:\},\:$ so $\:C_{ik}\:$ is contained in a class of 
$\:P^*(e_{ik})\:$ and, similarly, $\:\overline{C_{ik}}\:$ has such a
property. It remains to refer to their cardinalities again.
$\ \ \ \ \ \ \ \ \ \ \ \ \ \ \ \ \ \ \ \ \ \ \ \ \ \ \ \ \ \ \ \ \ \ \ \ \
\ \ \ \ \ \ \ \ \ \ \ \ \ \ \ \ \ \ \ \ \ \ \ \ \ \ \ \ \ \ \ \ \ \ \ \ \ 
\ \ \ \  \ \ \ \ \ \ \ \ \ \ \ \ \ \ \ \ \ \ \ \ \ \ \ \ \ \ \ \ \Box\ $

Now we able to obtain the multiplication table of the algebra $\:\CA.\:$

First of all we get
$$
e_{ik}e_{jk}\:=\:\frac{e_{ik}+e_{jk}}{2},\ \ \ 
e_{ik}e_{il}\:=\:\frac{e_{ik}+e_{il}}{2}
\eqno (5.11)
$$
from Lemma 5.12. Now if $\:i\neq j\:$ and $\:k\neq l\:$ we have $\:R(e_{ik},
e_{jl})\:\neq\:0\:$ from the same lemma which says that $\:e_{jl}\:\not\in\:
P(e_{ik})\:$ in this case. Then there is a q.a. $\:\{\:e_{ik},e_{jl}
\mid e_g,e_h\:\}\:$ with
$$
\{\:e_g,e_h\:\}\:\subset\:P(e_{ik})\:\bigcap\:P(e_{jl})\:=\:
\{\:e_{il},e_{kj}\:\}.
$$
This means that
$$
e_{ik}e_{jl}\:=\:\frac{e_{il}+e_{kj}}{2}\ \ \ (i\neq j,\ k\neq l).
\eqno (5.12).
$$

The algebra $\:\CA\:$ turns out to be an e.q.a. Hence $\:\CA\:$ is regular.
The Main Theorem is proved in the case under consideration.

B) {\bf The mixed case} $\:(\rho <n).\:$ In this case we can argue by 
induction on $\:n.\:$ Recall that the Main Theorem is true for $\:n\leq 5\:$
by Theorem 2.14.

Given $\:n\geq 6,\:$ we suppose that the theorem is true in all dimensions 
less 
than $\:n,\:$ in particular, for all proper coordinate subalgebras. All of
them are ultranormal together with $\:\CA.\:$ Therefore they are regular 
and then each one is either u.a. or e.M.a., or e.q.a. (Theorem 3.1).
As a result, all proper coordinate subalgebras are grounded. By Corollary 
4.14 $\:\CA\:$ is also grounded, i.e. $\:\CA\:$ is the offspring 
subalgebra of the set of its basis idempotents, say $\:\{\:
e_i\:\}_1^{\rho}.\:$ 

{\bf Lemma 5.13} {\em For every basis vector $\:e_j\:$ with $\:j >\rho\:$
there exists a unique pair $\:\{\:e_{i_j},e_{k_j}\:\}\:$ with $\:1\leq i_j
<k_j\leq\rho\:$ such that {\rm Lin}$\{\:e_{i_j},e_{k_j},e_j\:\}\:$ is an 
e.M.a., $\:\{\:e_{i_j},e_{k_j}\mid e_j\:\}.\:$ }

{\bf Proof.} Since $\:\CA\:$ is the offspring subalgebra of the family of 
idempotents $\:e_1,\:...,\:e_{\rho},\:$ there exists a pair 
$\:e_{i_j},e_{k_j}\ (1\leq i_j
<k_j\leq\rho)\:$ such that $\:e_j\:\in\:{\rm supp}(e_{i_j}e_{k_j}).\:$
The offspring subalgebra of this pair is neither unit (2-dimensional)
nor quadrille (because $\:e_j\:$ is not an idempotent). By Theorem 4.11 
it is an e.M.a. By (3.11)
$$
\{\:e_{i_j},e_{k_j}\:\}\:=\:{\rm supp}(e_j^2)\setminus\{\:e_j\:\},
$$
therefore the pair $\:\{\:e_{i_j},e_{k_j}\:\}\:$ is unique.
$\ \ \ \ \ \ \ \ \ \ \ \ \ \ \ \ \ \ \ \ \ \ \ \ \ \ \ \ \ \ \ \ \ \ \ \ 
\ \ \ \ \ \ \  \ \ \ \ \ \ \ \ \ \ \ \ \ \ \ \ \ \ \ \ \ \ \
\ \ \ \ \ \ \Box\ $

We will say that $\:e_j\:$ is the {\em offspring} of the {\em marked}
pair $\:\{\:e_{i_j},e_{k_j}\:\}.\:$ We also set $\:e_{ik}\:=\:e_ie_k\ \ 
(1\leq i\leq k\leq n),\:$ so that $\:e_{ii}\:=\:e_i,\ \ 1\leq i\leq\rho ,\:$
and Lin$\:\{\:e_{i_j},e_{k_j},e_{i_jk_j}\:\}\:=\:{\rm Lin}\:\{\:
e_{i_j},e_{k_j},e_j\:\}.\:$

{\bf Corollary 5.14.} {\em The algebra $\:\CA\:$ is nuclear.}

For any nonmarked pair $\:\{\:e_i,e_k\:\}\:$ of the basis idempotents the
offspring subalgebra is either u.a. or q.a.

{\bf Lemma 5.15.} {\em For every basis idempotent $\:e_i\:$ there are exactly
$\:\rho -m\:$ idempotents $\:e_k\ (k\neq i)\:$ such that the offspring 
subalgebra of the pair $\:\{\:e_i,e_k\:\}\:$ is q.a.} 

Thus, this number is the same for all $\:e_i.\:$

{\bf Proof.} Let $\:i=1\:$ for definiteness and let the offspring subalgebra
$\:L_k\:$ of the pair $\:\{\:e_1,e_k\:\}\:$ be an e.M.a., i.e.
$$
2e_1e_k\:=\:\ga_ke_1\:+\:\gb_ke_k\:+\:\gamma_ke_{\rho +k},\ \ \ 
\gamma_k >0.
$$
which means that $\:Be_k\:=\:\gb_ke_k\:+\:\gamma_ke_{\rho +k}\:$ where 
$\:B\:$ is the usual projection associated wich $\:e_1.\:$ If now $\:L_k\:$
is the u.a. then $\:Be_k=e_k.\:$ If, finally, $\:L_k\:$ is a q.a. then 
$\:Be_k\:=\: e_j+e_l\:$ where $\:e_j\:$ and $\:e_l\:$ correspond to the 
unit $\:L_j\:$ and $\:L_l.\:$ We see that
rk$\:\{\:Be_k\:\}_2^{\rho}\ {\rm equals}\:$ the total number of e.M.a. and 
u.a. On the other hand, this rank equals $\:m-1\:$ by Corollary 5.2. 
Therefore the number of q.a. coincides with $\:\rho -m.\:$
$\ \ \ \ \ \ \ \ \ \ \ \ \ \ \ \ \ \ \ \ \ \ \ \ \ \ \ \ \ \ \ \ \ \
\ \ \ \ \ \ \ \ \ \ \ \ \ \ \ \ \ \ \ \ \ \ \ \ \Box\ $

Lemma 5.15 is a quite preliminary fact because of

{\bf Lemma 5.16.} {\em There is no q.a. among the coordinate subalgebras of
$\:\CA.\:$ }

{\bf Proof.} Let we have the e.M.a. $\:\{\:e_1,e_2\mid e_{\rho +1}\:\}\:$
jointly with a q.a. By Lemma 5.15 $\:e_1\:$ is involved in a q.a., say,
$\:\{\:e_1,e_3\mid e_4,e_5\:\}.\:$ The offspring set of the triple
$\:\{\:e_1,e_2,e_3\:\}$ is 
$\{\:e_j\:\}_1^5\:\bigcup\:\{\:e_{\rho +1}\}\:\bigcup\:{\rm supp}(e_2e_3).\:$
 The corresponding offspring subalgebra is $\:\CA\:$ because of the
regularity of all proper coordinate subalgebras. 

Let the offspring subalgebra of the pair $\:\{\:e_2,e_3\:\}\:$ be u.a. or
e.M.a., so that there is no new idempotents in supp$(e_2e_3).\:$ Then
$\:\rho =5,\ e_{\rho +1}=e_6\ {\rm and}\ \{\:e_1,e_3\mid e_4,e_5\:\}\:$
is the unique q.a. containing $\:e_1.\:$ By Lemma 5.15 $\:\rho -m=1\:$
(so that $\:m=4)\:$ and then $\:e_2\:$ must be also involved in a q.a.
$\:\{\:e_2,e_i\mid e_k,e_j\:\}\:$ where $\:i=4\:$ or $\:i=5\:$ since the
offspring subalgebras of $\:\{\:e_1,e_2\:\}\:$ and $\:\{\:e_2,e_3\:\}\:$
are not q.a. We get at least two q.a. containing $\:e_4\ ({\rm or}\ e_5)\:$
in contradiction to Lemma 5.15.

Suppose that the offspring subalgebra of $\:\{\:e_2,e_3\:\}\:$ is a q.a.,
say, $\:\{\:e_2,e_3\mid e_i,e_k\:\},\:$ so that $\:\CA\:=\:{\rm Lin}(
\{\:e_j\:\}_1^5\:\bigcup\:\{\:e_i,e_k,e_{\rho +1}\:\} ).\:$
Now $\:e_3\:$ is involved in two q.a. but for $\:e_1\:$ there is no more 
q.a. than $\:\{\:e_1,e_3\mid e_4,e_5\:\}.\:$
Indeed, such a q.a. must be $\:\{\:e_1,e_i\mid *,*\:\}\:$ (up to 
transposition $\:e_i\leftrightarrow e_k).\:$ Both of the omitted members must
satisfy $\:R(e_1,*)=0\ {\rm and}\ R(e_i,*)=0.\:$ They must be $\:e_4\ 
{\rm or}\ e_5\:$ but $\:\{\:e_1,e_i\mid e_4,e_5\:\}\:$ contradicts the
pre-existence of $\:\{\:e_1,e_3\mid e_4,e_5\:\}.\:$
$\ \ \ \ \ \ \ \ \ \ \ \ \ \ \ \ \ \ \ \ \ \ \ \ \Box\ $

Lemmas 5.15 and 5.16 immediately imply

{\bf Corollary 5.17.} $\:\rho =m.\:$

Thus, the set of the basis idempotents is $\:\{\:e_i\:\}_1^m.\:$

{\bf Lemma 5.18.} {\em If $\:m\leq 4\:$ then $\:\CA\:$ is regular. }

By Corollary 5.14 and Theorem 2.11 $\:\CA\:$ is regular for types 
$\:(m,\gd)\:$ such that $\:m\leq 3\ {\rm or}\ \gd\leq 1,\ {\rm or}\ 
\gd\geq\frac{1}{2}(m-1)(m-2)+1.\:$ Thus, we can assume $\:m\geq 4\ {\rm or}
\ 2\leq\gd\leq\frac{1}{2}(m-1)(m-2),\:$ so that
$$
6\:\leq\:n\:\leq\:\frac{m(m-1)}{2}\:+\:1.
\eqno (5.13)
$$

{\bf Proof.} By (5.13) $\:n=6\:$ or $\:n=7.\:$ Let $\:n=6,\:$ so that type of
$\:\CA\:$ is (4,2). By Corollary 5.17 the basis idempotents are 
$\:e_1,e_2,e_3,e_4\:$ and there are exactly two marked pairs,
$\:\{\:e_{i_1},e_{k_1}\:\}\ {\rm and}\ \{\:e_{i_2},e_{i_{k_2}}\:\} ;\ 
e_5\ {\rm and}\ e_6\:$ are respectively their offsprings. Suppose that those
pairs do intersect, say, they are $\:\{\:e_1,e_2\:\}\ {\rm and}\ \{\:e_1,
e_3\:\} .\:$ By Lemma 5.16 the offspring subalgebras of all nonmarked pairs
$\:\{\:e_i,e_k\:\}\:$ are u.a. The offspring subalgebra of the triple
$\:\{\:e_1,e_2,e_3\:\}\:$ is
$$
L\:=\:{\rm Lin}\{\:e_1,e_2,e_3,e_5,e_6\:\}\:=\:{\rm Lin}
\{\:e_1,e_2,e_3,e_{12},e_{13}\:\}
$$
It is regular since dim$L=5.\:$ Now $\:\CA\:=\:L[e_4]\:$ is regular by 
Proposition 2.10, part 1.

If 
$\:\{\:e_{i_1},e_{k_1}\:\}\:\bigcap\:\{\:e_{i_2},e_{k_2}\:\}\:=\:\emptyset\:$ 
one can assume that those pairs are $\:\{\:e_1,e_2\:\}\ {\rm and}\ 
\{\:e_3,e_4\:\}.\:$ As before, we have the u.a. 
$\:\{\:e_i,e_k\mid\emptyset\:\}\:$ with $\:i=1,2\ {\rm and}\ k=3,4.\:$ The
offspring subalgebra $\:L\:=\:{\rm Lin}\{\:e_1,e_2,e_5\:\}\:$ is regular 
being e.M.a., and $\:\CA\:=\:L[e_3,e_4]\:$ is regular by Proposition 2.10,
part 2. (Condition (2.35) in the form $\:R(e_1,e_2)e_3=0\:$ is fulfilled 
because the offspring subalgebra Lin$\{\:e_1,e_2,e_3,e_5\:\}\:$ is 
regular being of dimension 4.)

{\em On this stage the Main Theorem is proved for $\:n\leq 6.\:$ } 

Henceforth $\:n=7,\:$ so that $\:\CA\:$is of type (4,3). Then there are 
exactly three marked pairs, $\:\{\:e_{i_1},e_{k_1}\:\}\ \ 
\{\:e_{i_2},e_{k_2}\:\},\ \ \{\:e_{i_3},e_{k_3}\:\};\:$ their offsprings are
$\:e_5,e_6,e_7\:$ respectively. We have the u.a. $\:\{\:e_i,e_k\mid
\emptyset\:\}\:$ for all nonmarked pairs again.

Suppose the intersection of all marked pairs is not empty. Then they are
$\:\{\:e_1,e_2\:\},\ \ \{\:e_1,e_3\:\}\:$ and $\:\{\:e_1,e_4\:\}\:$ (up
to enumeration). Every triple $\:\{\:e_1,e_{ik},e_{i^{'}k^{'}}\:\}\:$ 
belongs to 
an offspring subalgebra of dimension$\:\leq 5\:$ which is regular a 
fortiori. Hence, $\:R(e_{ik},e_{i^{'}k^{'}})e_1\:=\:0\:$ for all $\:i,k,
i^{'},k^{'}\:$ (some of them may coincide). By Corollary 2.9 $\:\CA\:$ is 
regular.

Let the intersection of all marked pairs is empty. However, since all 
$\:e_{i_j}\:$ and $\:e_{k_j}\:$ are from $\:\{\:e_l\:\}_1^4,\:$ there 
exists a couple of the pairs with nonempty intersection, say 
$\:\{\:e_1,e_2\:\}\:$ and $\:\{\:e_1,e_3\:\}.\:$ Suppose that the third pair 
does
not include $\:e_4.\:$ Then this is $\:\{\:e_2,e_3\:\}.\:$ The offspring 
subalgebra $\:M\:=\:{\rm Lin}\{\:e_1,e_2,e_3,e_5,e_6,e_7\:\}\:$ is 
regular because of dim$L=6.\:$ Since $\:\CA\:=\:\bM[e_4],\:$ this is 
regular by Proposition 2.10, part 1.

Suppose that $\:e_4\:$ is included in the third pair. Then those pairs are
$\:\{\:e_1,e_2\:\},\:$ $\\$
$\:\{\:e_1,e_3\:\}\ {\rm and}\ \{\:e_3,e_4\:\}\:$ (up to the
transposition $\:e_2\leftrightarrow e_3).\:$ We are going to use 
Corollary 2.9 again. In this context we only must consider those triples 
$\:\{\:e_1,e_{ik},e_{i^{'}k^{'}}\:\}\:$ wich are not located in a proper 
offspring
subalgebra because all of these subalgebras are of dimension$\:\leq\:6.\:$
Such "badly located" triples appear iff $\:e_2,e_3\:$ and $\:e_4\:$ are all 
among $\:e_i,e_k,e_{i^{'}},e_{k^{'}}.\:$ Because of the symmetry between 
$\:e_{ik}\:$ and $\:e_{i^{'}k^{'}}\:$ one can assume that $\:\{\:i,k\:\}\:$
is a lexicographic predecessor of $\:\{\:i^{'},k^{'}\:\}.\:$ The complete 
list of couples $\:\{\:e_{ik},e_{i^{'}k^{'}}\:\}\:$ under our conditions is
the following:
$$
\{ e_{12},e_{34}\},\ \{ e_{13},e_{24}\},\ \{ e_{14},e_{23}\},\ \{ e_{22},
e_{34}\},\ \{ e_{23},e_{24}\},\ \{ e_{23},e_{34}\},\ \{ e_{24},e_{33}\},\ 
\{ e_{24},e_{34}\}.
$$
Here $\:e_{jl}\:=\:\frac{1}{2}(e_j +e_l)\:$ (i.e. $\:R(e_j,e_l)=0)\:$ except
for the marked ones, i.e. $\:e_{12},e_{13}\:$ and $\:e_{34}.\:$ In 
particular,$\:e_{14}\:=\:\frac{1}{2}(e_1+e_2)\:$ and 
$\:e_{23}\:=\:\frac{1}{2}(e_2+e_3).\:$ Hence,
$$
R(e_{14},e_{23})\:=\:\frac{1}{4}\{\:R(e_1,e_2)\:+\:R(e_1,e_3)\:+\:R(e_3,e_4)
\}
\eqno (5.14)
$$
because of $\:R(e_2,e_4)=0.\:$ Since every of triples 
$\:\{\:e_1,e_1,e_2\:\},\ 
\{\:e_1,e_1,e_3\:\}\ {\rm and}\ \{\:e_1,e_3,e_4\:\}\:$ is well located, 
(5.14) yields 
$$
R(e_{14},e_{23})e_1\:=\:0.
\eqno (5.15)
$$
Similarly, $\:R(e_{23},e_{24})\:=\:\frac{1}{4}R(e_3,e_4)\:=\:R(e_{24},e_{33})
\:$ hence,
$$
R(e_{23},e_{24})e_1\:=\:0,\ \ R(e_{24},e_{33})e_1\:=\:0.
\eqno (5.16)
$$
Now $\:R(e_{22},e_{34})\:=\:R(e_2,e_{34})\:=\:-\frac{1}{2}R(e_3,e_4)\:$ by
(2.25). Hence,
$$
R(e_{22},e_{34})e_1\:=\:0.
\eqno (5.17)
$$
Since $\:R(e_{23},e_{34})\:=\:\frac{1}{2}\{ 
R(e_2,e_{34})\:+\:R(e_3,e_{34})\},\:$ we obtain from (5.17)
$$
R(e_{23},e_{34})e_1\:=\:\frac{1}{2}R(e_3,e_{34})e_1\:=\:0
\eqno (5.18)
$$
because the triple $\:\{\:e_1,e_3,e_{34}\:\}\:$ is well located. Similarly,
$$
R(e_{24},e_{34})e_1\:=\:0.
\eqno (5.19)
$$
Since $\:R(e_{13},e_{24})\:=\:\frac{1}{2}\{ R(e_{13},e_2)\:+\:R(e_{13},e_4)\}
\:$ and the triples $\:\{\:e_1,e_{13},e_2\:\}\:$ and 
$\:\{\:e_1,e_{13},e_4\:\}\:$ are well located, we obtain
$$
R(e_{13},e_{24})e_1\:=\:0
\eqno (5.20)
$$
Finally, $\:R(e_{12},e_{34})\:=\:-\{ R(e_{13},e_{24})\:+\:R(e_{14},e_{23})\}
\:$ by (2.20) and then
$$
R(e_{12},e_{34})e_1\:=\:0
\eqno (5.21)
$$
by (5.15) and (5.20).
$\ \ \ \ \ \ \ \ \ \ \ \ \ \ \ \ \ \ \ \ \ \ \ \ \ \ \ \ \ \ \ \ \ \ \ \ 
\ \ \ \ \ \ \ \ \ \ \ \ \  \ \ \ \ \ \ \ \ \ \ \ \ \ \ \ \ \ \ \ \ \ \ \ 
\ \ \ \ \ \ \ \ \ \ \ \ \ \ \ \ \Box\ $

Now we are able to finish the proof of the Main Theorem. Actually, we
are going to prove that $\:\CA\:$ is an e.M.a., which means that 
(3.1)-(3.5) is valid for the set 
$\:\{\:e_i\:\}_1^m\:$ of the basis idempotents. We already have (3.1) 
with marked pairs $\:\{\:e_{i_j},e_{k_j}\:\}\:$ coming from Lemma 5.13
(where $\:\rho =m\:$ by Corollary 5.17). For nonmarked pairs $\:\{\:e_i,e_k\:
\}\:$ we have (3.2) because of Theorem 4.11 and Lemma 5.16.

In order to get (3.3) and (3.4) we consider the offspring subalgebra
$\:L_{ij}\:$ of the set $\:\{\:e_i,e_{i_j},e_{k_j}\:\}.\:$
Since dim$L\leq 6\:$ (the equality is attained if both of pairs
$\:\{\:e_i,e_{i_j}\:\}\:$ and $\:\{\:e_i,e_{k_j}\:\}\:$ are marked ), $\:L\:$
is regular. Therefore we have (3.3) and (3.4), moreover, $\:c_j\:$ 
is independent of $\:i\:$ by virtue of (3.4) and also $\:\overline{c_j}\:=\:
1-c_j.\:$ (Recall that $\:\gamma_j >0.)\:$

It remains to get (3.5). For this goal we consider the offspring subalgebra
$\:M_{jl}\:$ of the set $\:\{\:e_{i_j},e_{k_j},e_{i_l},e_{k_l}\:\}.\:$
The offspring subalgebras of all six pairs $\:\{\:e_{i_j},e_{k_j}\:\},\ 
\{\:e_{i_j},e_{i_l}\:\},\ ...,\ \{\:e_{i_l},e_{k_l}\:\}\:$ are e.M.a. or 
u.a. Therefore the idempotents in the canonical basis of $\:M_{jl}\:$ are 
only
$\:e_{i_j},e_{k_j},e_{i_l}\:$ and $\:e_{k_l}\:$ i.e. $\:\rho\leq 4\:$ for 
$\:M_{jl}.\:$ By Corollary 5.17 $\:m\leq 4.\:$ By Lemma 5.18 $\:M_{jl}\:$
is regular. Hence, (3.5) is also valid.

\baselineskip  18pt
\medskip
\bigskip
\ct{\bf References}
\bigskip
\bigskip

 1. S.N. Bernstein. Mathematical problems in modern biology. {\em
Science in the Ukraine } {\bf 1} (1922), 14-19 (in Russian). 

 2. S.N. Bernstein. Demonstration mathematique de la loi d'heredite de Mendel.
{\em C.r.Acad. Sci. Paris} {\bf 177} (1923), 528-531.

 3. S.N. Bernstein, Principe de stationarite et generalisation de la loi
de Mendel. {\em C.r.Acad. Sci. Paris} {\bf 177} (1923), 581-584.

 4. S.N. Bernstein. Solution of a mathematical problem related to the theory
of inheritance. {\em Uch. Zap. n.-i. kaf. Ukrainy}, {\bf 1} (1924),
83-115 (in Russian).

 5. G.H. Hardy. Mendelian proportions in a mixed population. {\em Science}
{\bf 28}, 706 (1908), 49-50.

 6. P. Holgate. Genetic algebras satisfying Bernstein stationarity principle.
{\em J.London Math. Soc.} {\bf 9} (1975), 613-624.

 7. Yu.I. Lyubich. Basic concepts and theorems of evolutionary genetics for
free populations. {\em Russian Math. Surveys} {\bf 26}, 5 (1971), 51-123.

 8. Yu.I. Lyubich. On the mathematical theory of heredity. {\em Soviet Math.
Dokl.} {\bf 14}, 2 (1973), 579-581.

 9. Yu.I. Lyubich. About a theorem of S.N. Bernstein. {\em Siberian Math. J.}
{\bf 14}, 3 (1973), 474.

10. Yu.I. Lyubich. A class of quadratic maps. {\em Teor.
Func.,Func. Anal. and Appl.} {\bf 21} (1974), 36-42 (in Russian).

11. Yu.I. Lyubich. Two-level Bernstein populasions. {\em Math. of the USSR
Sbornik} {\bf 24}, 1 (1974), 593-615.

12. Yu.I. Lyubich. Proper Bernstein populations. {\em Probl. Peredachi
Inform. } {\bf 13}, 3 (1977), 91-100. Transl. in {\em Probl. 
Inform.Transmiss.} (1978), Jan., 228-235.

13. Yu.I. Lyubich. Bernstein algebras. {\em Uspekhi Mat. Nauk}
{\bf 32}, 6 (1977), 261-263 (in Russian).

14. Yu.I. Lyubich. Stochastic algebras and their applications to mathematical
genetics, In {\em Mathematical methods in biology.}
Naukova Dumka, Kiev, (1977), 119-131 (in Russian).

15. Yu.I. Lyubich. Algebraic proof of S.N.Bernstein's theorem on two pure
types. {\em Vestnik Kharkov Univ., Ser. Mat. and Mech.}
{\bf 177}, 44 (1979), 86-94 (in Russian).

16. Yu.I. Lyubich. A topological approach to a problem in mathematical 
genetics. {\em Russian Math. Surveys} {\bf 34} (1979), 6, 60-66.

17. Yu.I. Lyubich. {\em Mathematical Structures in Population Genetics.}
Naukova Dumka, Kiev, 1983. English transl. in Springer, 1992.

18. W. Weinbreg. Uber den Nachweis der Verebung beim Menschen.
{\em  Jahresber. Ver. vaterl. Naturk. in Wurtemb.} {\bf 64} (1908), 368-372.
\\
\\


\end{document}